\documentclass[a4paper,12pt]{article}
\usepackage{amssymb,amsmath}
\renewcommand{\mod}[1]{\hspace{-2.9mm}\pmod{#1}}

\newcommand{\li}{{\rm  Li}}

\newcommand{\twosum}[2]{\sum_{\begin{array}{c}{\scriptstyle #1}\\
{\scriptstyle #2} \end{array}}}
\newcommand{\twoprod}[2]{\prod_{\begin{array}{c}{\scriptstyle #1}\\
{\scriptstyle #2} \end{array}}}
\newcommand{\twoar}[2]{\big({\begin{array}{c} { #1}\\
{ #2} \end{array}}\big)}

\newcommand{\twos}[4]{\Big\{{\begin{array}{cc} { #1} & { #2}\\
{ #3} & { #4} \end{array}}}
\newcommand{\bg}[1]{\big( { #1} \big )}
\newcommand{\Bg}[1]{\Big( { #1} \Big )}

\newcommand{\R}{\mathbb{R}}

\newcommand{\N}{\mathbb{N}}
\newcommand{\op}{\mathbb{P}}
\renewcommand{\b}[1]{{\bf #1}}

\newcommand{\cl}[1]{{\cal #1}}

\newcommand{\ep}{\varepsilon}
\newcommand{\Z}{\mathbb{Z}}

\newtheorem{theorem}{Theorem}[section]
\numberwithin{equation}{section}
\newtheorem{lemma}{Lemma}[section]
\newtheorem{cor}{Corollary}[section]

\begin{document}

\title{Lectures on sieves}
\author{D.R. Heath-Brown\\Mathematical Institute, Oxford}
\date{}
\maketitle

\renewcommand{\thesection}{}
\section{Preface}
\renewcommand{\thesection}{1}

These are notes of a series of lectures on sieves, presented during
the Special Activity in Analytic Number Theory, at the Max-Planck
Institute for Mathematics in Bonn, during the period January--June
2002.  The notes were taken by Boris Moroz, and it is a pleasure to
thank him for writing them up in the current form.  In addition,
thanks are due to the institute for its hospitality and financial
support, and also to the American Institute of Mathematics, where these
notes were put into final form.

Being lecture notes they are not intended to have the formal style of
a textbook.  Nor is there any significant claim to novelty.  In particular the
material in \S\S 3 \& 4 owes much to the book of Halberstam and
Richert \cite{HR}.  In addition to the latter work the reader may find it
helpful also to consult Greaves' exposition \cite{GreavesBook}.

\section{Introduction}
Here are some typical questions which can be expressed as sieve
problems:

(i) Is every even integer $n\ge 2$  a sum of two
primes? ({\it Goldbach's conjecture});

(ii) Are there infinitely many pairs of primes $(p, q)$ with
$q=p+2$? ({\it the twin primes conjecture});

(iii) Are there arbitrary long arithmetic progressions consisting only of
primes ? ;

(iv) Are there infinitely many primes of the form $n^2 + 1$ with
$n\in \N$ ? ;

(v) Is it true that, for every $n\in \N$, there is a prime $p$
in the range $n^2<p<(n+1)^2$ ? ;

(vi) For every $\ep>0$ is there an integer $N(\ep)$ such that the
interval $[N,N+N^{\ep}]$ contains a square-free number, as soon as
$N\ge N(\ep)$? ;

(vii) Is it true that, for some
$\gamma\leq 1/4$ and every sufficiently large $x$,  there is a
pair of natural numbers $(m, n)$ such that $x\leq m^{2}+n^{2} <
x+x^{\gamma}$ ?.

{\bf Exercise.} Prove that
$$\{(m, n): {m, n}\in \N,\;x\leq m^{2}+n^{2} < x+8x^{1/4}\}
\neq \emptyset$$ for all sufficiently large $x$.
\bigskip

Let $\cl{A}$ be a finite subset of $\N$, and $\cl{P}\subset \op$,
where $\op$ denotes the set of all primes. 
For any positive real $z$, set
$$S(\cl{A}, \cl{P}; z) :=\#\;\{n\in \cl{A}:\; p|n \Rightarrow p\geq
z\; \text{for}\; p\in \cl{P}\}.$$

{\bf Example 1.} Let $\cl{A} =\{n\in \N:\; n\leq
x\}$ and $\cl{P}=\op$. Then
$$S(\cl{A}, \cl{P}; z) :=\#\;\{n\in \N:\; n\leq x,\; p|n\Rightarrow
p\geq z\; \text{for}\; p\in \op\},$$
so that $S(\cl{A}, \cl{P}; z) = \pi (x) - \pi (z)$ for
$\sqrt{x}<z\leq x$, where, as usual,
$$\pi (x):=\#\;\{p\in \op:\; p\leq x\}.$$

{\bf Example 2.} Let $\cl{A} = \{n(2N-n): n\in \N,\; 2\leq n\leq
2N-2\}$ and $\cl{P}=\op.$ Then
\begin{eqnarray*}
\lefteqn{S(\cl{A}, \cl{P}; \sqrt{2N})}\\
& =&\#\;\{ p\in \op:\; 2N-p\in \op,\;\sqrt{2N}\leq p, 2N-p\leq 2N-2\}.
\end{eqnarray*}
This relates to Goldbach's problem.

{\bf Example 3.} Let
$$\cl{A} =\{n\in \N:\; n<x\},\;\;\;
\cl{P}=\{p\in \op:\; p\equiv 3\mod{4}\}.$$ Then
$$S(\cl{A}, \cl{P}; x) :=
\#\;\{n\in \N:\; n< x,\; p\nmid n\; \text{for}\;
p\in \op,\;p\equiv 3\mod{4}\},$$ so that
$$S(\cl{A}, \cl{P}; x) = \#\;\{n\in \N:\; n<x,\; n=l^{2}+m^{2}\; 
\text{with}\;l, m \in \N,\; \text{h.c.f.}(l, m) = 1\}.$$
We can therefore detect sums of two squares.
\bigskip

Let now $\Pi := \Pi (\cl{P}, z): = \prod\limits_{p<z,\, p\in \cl{P}}
p$, and let $\cl{A}_{d} :=\{n\in \N:\; nd\in \cl{A}\}$. It
follows that
\begin{equation}
S(\cl{A}, \cl{P}; z) = \twosum{n\in \cl{A}}{(n, \Pi)=1} 1 =
\twosum{n\in \cl{A}}{d|(n, \Pi)} \mu (d) = \sum_{d|\Pi (\cl{P},
z)} \mu (d) \#\;\cl{A}_{d}.
\end{equation}
{\bf Main assumption.}  {\it Suppose that}
\begin{equation}
\#\;\cl{A}_{d}=\frac{\omega(d)}{d}X + R_{d},
\end{equation}
{\it with an absolutely multiplicative $\omega(d)$, satisfying the
conditions $\omega(d)\geq 0$ for $d\in \N$ and
$\omega(p)=0$ for $p\in \op\setminus \cl{P}$.}

In particular, $\#\;\cl{A}=X + R_{1}$. We think of the remainder
$R_{d}$ term as being small compared to the main term
$\frac{\omega(d)}{d}X$.

In the notation of Example 2, one may write
\begin{eqnarray*}
\#\;\cl{A}_{d} &=& \#\{n(2N-n): n\in \N,\; 2\leq n\leq 2N-2,\;
d|n(2N-n)\}\\
&=&\twosum{m\mod{d}}{d|m(2N-m)}\#\{n:2\le n\le 2N-2,\;n\equiv m\mod{d}\}\\
& =&\twosum{m\mod{d}}{d|m(2N-m)}
(\frac{2N}{d} + r(d,m)),
\end{eqnarray*}
where $r(d,m)\ll 1$.  Thus if $\omega(d)=\#\;\{m\mod{d}:d|m(2N-m)\}$, then
\[\#\;\cl{A}_{d}=\frac{\omega(d)}{d}X + R_{d}\]
with $X=2N$ and so, for Example 2, we have
\[R_d=\twosum{m\mod{d}}{d|m(2N-m)}r(d,m)\ll\omega(d).\]
\bigskip

In general it follows from (1.1) and (1.2) that
\begin{equation}
S(\cl{A}, \cl{P}; z) = X \sum_{d|\Pi (\cl{P}, z)} \frac{\mu(d)
\omega(d)}{d} + \sum_{d|\Pi (\cl{P}, z)} \mu(d) R_{d}.
\end{equation}
Writing $W(z; \omega)=\prod\limits_{p<z,\, p\in \cl{P}} (1 -
\frac{\omega(p)}{p})$ we deduce the following from (1.3).
\begin{cor}
We have
\[
|S(\cl{A}, \cl{P}; z) -  X W(z; \omega)|\leq \sum_{d|\Pi (\cl{P},
z)} |R_{d}|.
\]
\end{cor}
This corollary is known as the {\it Sieve of Eratosthenes$\,$--Legendre}.

\begin{lemma}
{\bf The Mertens Formula.}
Let $V(z):= \prod\limits_{p<z,\, p\in \op}(1-1/ p)$. One has
\begin{equation}
V(z) = \frac{e^{-\gamma}}{\log z} + O(\frac{1}{(\log z)^2}).
\end{equation}
\end{lemma}
{\bf Proof.} See, for instance, Prachar \cite[pp. 80-81.]{Prach}

{\bf Example 1} (continued). We have
$$\cl{A}_{d} =\{n\in \N:\; nd\leq x\},\;\;\;\#\cl{A}_{d} = [\frac{x}{d}] =
\frac{X}{d} + R_{d}$$
and $\cl{P}=\op,$ so that we may take $X = x,\;\omega(d) =
1,\; |R_{d}|\leq 1$. Therefore it follows from Corollary 1.1 and
Lemma 1.1 that
$$S(\cl{A}, \cl{P}; z) = x \frac{e^{-\gamma}}{\log z}\{1 + O(1/(\log z))\} +
O(2^{\pi(z)}),$$ since
$$
\sum_{d|\Pi (\cl{P}, z)} |R_{d}|\leq \sum_{d|\Pi (\cl{P}, z)} 1
= 2^{\pi(z)}.$$ Thus
$$S(\cl{A}, \cl{P}; z) \sim x \frac{e^{-\gamma}}{\log z}$$
as soon as $z\leq \log x$, since
\begin{eqnarray*}
2^{\pi(z)}& = &e^{\pi(z)\log 2}\\
& =& \exp (\frac{z}{\log z}\{1 + O(1/\log z)\}\log 2)\\
&\leq &\exp (\frac{\log x}{ \log\log x}\{1 + O(1/\log\log x)\}\log 2)\\
&\le & \exp(\frac{1}{2} \log x)\\
& = &x^{1/2}\\
&=&o(x/ \log z)
\end{eqnarray*}
for $z\leq \log x,\;x\rightarrow \infty$.
\begin{cor} For $z\le \log x$ we have
$$\#\;\{n\in \N:\; n\leq x,\; p|n \Rightarrow p\geq z\;
\text{for}\; p\in \op\}\sim \frac{ e^{-\gamma}x}{\log z}.$$
\end{cor}

We would like to extend the admissible range of $z$ in such a result.
The proof of Corollary 1.2 uses the fact that
\begin{equation*}
\sum_{d|n,\;d|\Pi (\cl{P}, z)}\mu(d)=\left\{\begin{array}{cc} 1, &
(n, \Pi (\cl{P}, z)) = 1 ,\\ 0, &  (n, \Pi (\cl{P}, z)) > 1.
\end{array}\right.
\end{equation*}
However the sum here is over an uncomfortably large range.  We
therefore replace the above equality with two inequalities, and
encounter the following.

\noindent {\bf Sieve problem.} Find two real-valued functions
$\mu^{+}(d)$ and $\mu^-(d)$, of suitably small support, 
satisfying the conditions
\begin{equation}
\sum_{d|n,\;d|\Pi (\cl{P}, z)}\mu^{-}(d)\leq
\left\{\begin{array}{cc} 1, & (n, \Pi (\cl{P}, z)) = 1 ,\\ 0, &
(n, \Pi (\cl{P}, z)) > 1, \end{array}\right.
\end{equation}
and
\begin{equation}
\sum_{d|n,\;d|\Pi (\cl{P}, z)}\mu^{+}(d)\geq
\left\{\begin{array}{cc} 1, & (n, \Pi (\cl{P}, z)) = 1 ,\\ 0, &
(n, \Pi (\cl{P}, z)) > 1. \end{array}\right.
\end{equation}
It follows from (1.6) for example, that
\begin{eqnarray}
S(\cl{A}, \cl{P}; z) &=& \sum_{n\in \cl{A},\;(n, \Pi (\cl{P}, z))=1}
1\nonumber\\
&\leq&
\twosum{n\in \cl{A},\;d|n}{d|\Pi (\cl{P}, z)} \mu^{+}(d)\nonumber\\
& =&
\sum_{d|\Pi (\cl{P}, z)} \mu^{+}(d) \#\;\cl{A}_{d}\nonumber\\
&=&
X \sum_{d|\Pi (\cl{P}, z)} \frac{\mu^{+}(d) \omega(d)}{d} +
\sum_{d|\Pi (\cl{P}, z)} \mu^{+}(d) R_{d}.
\end{eqnarray}
Hence
\begin{equation}
S(\cl{A}, \cl{P}; z)\leq X \sum_{d|\Pi (\cl{P}, z)}
\frac{\mu^{+}(d) \omega(d)}{d} + \sum_{d|\Pi (\cl{P}, z)}
|\mu^{+}(d)| |R_{d}|.
\end{equation}
To minimise the right-hand side of (1.8), subject to the condition
(1.6), is in general a challenging unsolved problem.

{\bf Some achievements of sieve methods.}

(i) We have
$$\#\;\{(p, p^\prime): p, p^\prime\in \op,\; p + p^\prime = 2n\}\ll
\frac{\sigma(n)}{n}\frac{n}{(\log n)^2} $$ where $\sigma(n)$ is the
sum of divisors function.  This is conjectured to be best possible, up
to the value of the implied constant.

(ii)  We have (Chen \cite{Chen})
$$\#\;\{(p, p^\prime): p\in \op,\;p^{\prime}\in \op_2,\; p + p^\prime = 2n\}\gg
\frac{n}{(\log n)^2},$$ where $\op_2$ is the set of positive integers
which are either prime or a product of two primes.

(iii) We also have (Chen \cite{Chen})
$$\#\;\{p\in \op:\;p\le x,\; p + 2\in \op_2\}\gg \frac{x}{(\log x)^2}.$$

(iv) We have (Iwaniec \cite{Iwann2+1})
$$\#\;\{n\le x: n^2 + 1\in \op_2\}\gg
\frac{x}{\log x}.$$

(iv) We have (Heath-Brown \cite{HBP2AP})
$$\#\;\{n \equiv l\mod{k}:\;n\in \op_2,\;n\leq k^{2}\}\gg
\frac{k^2}{\phi(k)\log k}$$
for large enough $k$, if ${\rm h.c.f.}(l, k)=1$.

\renewcommand{\thesection}{2}
\section{Selberg's sieve}

To satisfy (1.6), let
\begin{equation}
\mu^{+}(d) = \sum_{d = [d_{1}, d_{2}]}
\lambda_{d_{1}}\lambda_{d_{2}}
\end{equation}
with $\lambda_{1}=1$ and $\lambda_{d}\in \R$.
Clearly,
$$\sum_{d|n,\;d|\Pi (\cl{P}, z)}\mu^{+}(d) =
\sum_{[d_{1}, d_{2}]|(n, \Pi (\cl{P}, z))}
\lambda_{d_{1}}\lambda_{d_{2}}$$
$$ = \sum_{d_{1},d_2|(n, \Pi (\cl{P},
z))} \lambda_{d_{1}}\lambda_{d_{2}} =
(\sum_{d|(n, \Pi (\cl{P}, z))} \lambda_{d})^2\geq 0.$$
Moreover, if $(n, \Pi (\cl{P}, z))=1$ then
$$\sum_{d|n,\;d|\Pi (\cl{P}, z)}\mu^{+}(d) = \mu^{+}(1)=\lambda_{1}^2 = 1.$$
Hence $\mu^{+}(d)$ satisfies (1.6).

We shall minimise the main term on the right-hand side of (1.8) in the class of
functions given by (2.1). Let 
$$S_{0}:=\sum_{d|\Pi (\cl{P}, z)}\frac{\omega(d)\mu^{+}(d)}{d}\;,\;\;\;
\tilde{R}:=\sum{|\mu^{+}(d)|\,|R_{d}|}\,,\;\;\;\xi:=\sqrt{y}\,.$$
If $\lambda_{d}=0$ for $d\geq\xi\,,$ then $\mu^{+}(d)=0$ for $d\geq{y}\,.$
Now, 
$$S_{0}=\sum_{[d_{1},\, d_{2}]|\Pi (\cl{P}, z)}
\frac{\omega([d_{1},d_{2}])}{[d_{1},d_{2}]}\lambda_{d_{1}}\lambda_{d_{2}}\;.$$
Suppose that $\lambda_{d}=0$ for $d\geq\xi$, then 
$$S_{0}\,=\,{\twosum{d_{1},d_2|\Pi}{d_{1},d_2<\xi}}
\hspace{-4mm}^{*}
{\frac{\omega(d_{1})\lambda_{d_{1}}}{d_1}
\frac{\omega(d_{2})\lambda_{d_{2}}}{d_{2}}}\,
\frac{(d_{1},d_{2})}{\omega((d_{1},d_{2}))}\;,$$
where $\sum^{\,*}$ omits the terms for which
$\omega(d_{1})\,\omega(d_{2})=0$. 

We now introduce the assumption
$$0\leq \omega(p)< p\quad \mbox{ for}\quad 
p \in \cl{P}\quad\quad (\Omega_{1}),$$
which we shall refer to in future merely as {\bf Condition} $(\Omega_1)$.
Then if $\mu(k)\omega(k)\not=0$ we find that
\begin{eqnarray*}
\sum_{l|k}\mu\big( \frac{k}{l} \big)\,\frac{l}{\omega(l)}&=&\sum_{l'|k}
\mu(l')\,
\frac{k/l'}{\omega(k/l')}\\
&=&\frac{k}{\omega(k)}\,
\sum_{l'|k} \mu(l')\,\frac{\omega(l')}{l'}\\
&=&\frac{k}{\omega(k)}\, \prod_{p|k}
\big(1\,-\,\frac{\omega(p)}{p}\big)\\
&>& 0.
\end{eqnarray*}
Thus we may define a non-negative multiplicative function $g(k)$ by
\[\frac{1}{g(k)}=\sum_{l|k}\mu
\big(\frac{k}{l}\big)\,\frac{l}{\omega(l)}.\] 
The M\"{o}bius inversion formula then shows that
$$\sum_{k|d} \frac{1}{g(k)}\,=\, \frac{d}{\omega(d)}$$ 
for $\mu(d)\omega(d)\not=0$.  Under condition $(\Omega_{1})\,,$ it
now follows that 
\begin{eqnarray*}
S_{0}&=&{\twosum{d_{1},d_2|\Pi}{d_{1},d_2<\xi}}
\hspace{-6mm}^{*}\,
\frac{\omega(d_{1})\lambda_{d_{1}}}{d_{1}} \,
\frac{\omega(d_{2})\lambda_{d_{2}}}{d_{2}}\,
\sum_{l|(d_{1},d_{2})} \frac{1}{g(l)}\\
&=&\twosum{l|\Pi(\cl{P},\,z)}{l<\xi} \frac{1}{g(l)}
\,\Big(\twosum{d|\Pi(\cl{P},\,z)}{l|d,\,d<\xi} 
\frac{\omega(d)\lambda_{d}}{d}\;\,\Big)^{2}\\
&=&\twosum{l|\Pi(\cl{P}\,z)}{l<\xi} \frac{1}{g(l)}\,y^{2}_{l}
\end{eqnarray*}
with
$$y_{l}\,:=\,\twosum{d|\Pi(\cl{P},\,z)}{l|d,\,d<\xi}
\frac{\omega(d)\lambda_{d}}{d}\;.$$
\begin{lemma}
We have
\begin{equation}
\twosum{l|\Pi,\,d|l}{l<\xi}\mu(l)\,y_{l}\,=
\,\frac{\omega(d)\lambda_{d}}{d}\,\mu(d)\;,
\end{equation}
if $d<\xi$ and $d|\Pi$.
\end{lemma}
{\bf Proof.} When $d|\Pi$ we have  
\begin{eqnarray*}
\twosum{l|\Pi,\,l<\xi}{d|l}\mu(l)y_{l}&=&
\twosum{l|\Pi,\,l<\xi}{d|l}\twosum{\delta|\Pi,\,\delta<\xi}{l|\delta}
\mu(l)\,\frac{\omega(\delta)\lambda_{\delta}}{\delta}\\
&=& \twosum{\delta|\Pi,\,\delta<\xi}{d|\delta}
\frac{\omega(\delta)\lambda_{\delta}}{\delta}\,\sum_{d|l,\;l|\delta}\mu(l)\\
&=&\twosum{\delta|\Pi,\,\delta<\xi}{d|\delta}
\frac{\omega(\delta)\lambda_{\delta}}{\delta}\,\sum_{md|\delta}\mu(md)\\
&=&\mu(d)\,\twosum{\delta<\xi,\,\delta|\Pi}{d|\delta}
\frac{\omega(\delta)\lambda_{\delta}}{\delta}\,
\sum_{m|\frac{\delta}{d}}\mu(m)\\
&=&\frac{\mu(d)\,\omega(d)\,\lambda_{d}}{d}\;,
\end{eqnarray*}
as claimed.  Here we have used the fact that $\Pi$ is
square-free, so that $m$ and $d$ are coprime for $md|\Pi$.

Since $\lambda_{1}=1$, it follows from (2.2) that
$$1=\sum_{l|\Pi(\cl{P},\,z),\,l<\xi}\mu(l)\,y_{l}\,=\,
\sum_{l|\Pi(\cl{P},\,z),\,l<\xi}\mu(l)\sqrt{g(l)}\,
\frac{y_{l}}{\sqrt{g(l)}}$$ and therefore
$$1\le\big\{\,\sum_{l|\Pi(\cl{P},\,z),\;l<\xi}\mu(l)^{2}g(l)\,\big\}\,
\big\{\sum_{l|\Pi(\cl{P},\,z),\;l<\xi}y^{2}_{l}\,g(l)^{-1}\big\}=G(\xi,z)S_0$$
by Cauchy's inequality, where
$$G(\xi,\,z):=\sum_{l|\Pi,\,l<\xi}\mu(l)^{2}g(l) =
\sum_{l|\Pi(\cl{P},\,z),\,l<\xi}g(l)\;.$$ 
Thus $S_{0}\,\geq G(\xi,\,z)^{-1}\;,$ and
$S_{0}=G(\xi,\,z)^{-1}$ if and only if there is a constant 
$c$ such that 
$$\frac{y_{l}}{\sqrt{g(l)}}=c\mu(l)\,\sqrt{g(l)}\;$$ 
for every $l$, this being the condition for equality in Cauchy's 
inequality. For the optimal values
$y_{l}=c\mu(l)\,g(l)$ 
Lemma 2.1 yields
$$1\,=\,\frac{\omega(1)\lambda_1}{1}\mu(1)\,=\,
\sum_{l|\Pi,\,l<\xi}\mu(l)y_{l}\,=\,c\sum_{l|\Pi,\,l<\xi}
\mu(l)^{2}g(l)\,=\,c\,G(\xi,\,z),$$
whence
$$ c=\frac{1}{G(\xi,\,z)}\;.$$
Thus $S_{0}=G(\xi,\,z)^{-1}$ providing that
$y_{l}\,=\,\mu(l)\,g(l)G(\xi,\,z)^{-1}$, in which case (2.2) produces
\begin{eqnarray*}
\lambda_{d}& =& \mu(d)\,\frac{d}{\omega(d)}
\sum_{l|\Pi,\,d|l,\,l<\xi}
\mu(l)\,y_{l}\\
&=&\frac{\mu(d)\,d}{\omega(d)\,G(\xi,\,z)}
\twosum{l|\Pi,\,d|l}{l<\xi}\mu(l)^{2}g(l)\\
&=&\frac{\mu(d)\,d\,g(d)}{\omega(d)\,G(\xi,\,d)}
\sum_{dj|\Pi(\cl{P},\,z),\,j<\frac{\xi}{d}}g(j)\;.
\end{eqnarray*}
On recalling that
$$g(d)\,=\,
\frac{\omega(d)}{d}\,\prod_{p|d}\big(1-\frac{\omega(p)}{p}\big)^{-1}\;,$$
one obtains the minimising condition 
\begin{equation}
\lambda_{d}\,=\,\mu(d)\prod_{p|d}\big(1-\frac{\omega(p)}{p}\big)^{-1}\,
\frac{G_{d}(\frac{\xi}{d},\,z)}{G(\xi,\,z)}
\quad \mbox{for} \quad d<\xi\;,
\end{equation}
with
$$G_{d}(\xi,\,z):=\sum_{dj|\Pi(\cl{P},z),\,j<\xi}g(j)\;,\quad 
G_{1}(\xi,\,z)=G(\xi,\,z)\;.$$
The choice of variables (2.3), under the assumption
\begin{equation}
\lambda_{d}=0 \quad \mbox{for} \quad d\geq \xi\;,
\end{equation}
turns (1.8) into the inequality
\begin{equation}
S(\cl{A},\cl{P};z)\leq
X\,G(\xi,\,z)^{-1}\,+\,\sum_{d|\Pi(\cl{P},z),\,d<y}
|\mu^{+}(d)|\,|R_{d}|,
\end{equation}
since $\mu^{+}(d)=0$ for $d\geq y$, in view of (2.4).

In order to bound $\lambda_d$ we will require the following result.
\begin{lemma}
Let $d|\Pi(\cl{P},z)\,.$ Then
\begin{equation}
G(\xi,\,z)\geq
G_{d}\big(\frac{\xi}{d},\,z\big)\,\prod_{p|d}
\big(1-\frac{\omega(p)}{p}\big)^{-1}\;.
\end{equation}
\end{lemma}
{\bf Proof.} Let $d|\Pi(\cl{P},z).$ Then
\begin{eqnarray*}
G(\xi,\,z)&=&\sum_{l|d}\twosum{m<\xi}{m|\Pi,\,(m,\,d)=l} g(m)\\
&=&
\,\twosum{l|d,\,lh|\Pi}{lh<\xi,\,(h,\,d/l)=(h,\,l)=1}
g(lh)\\
&=&\sum_{l|d}g(l)\twosum{h<\xi/l,\,lh|\Pi}{(h,\,d/l)=(h,\,l)=1}g(h)\\
&=&\sum_{l|d} g(l)\sum_{h<\xi/l,\,dh|\Pi}g(h)\\
&\geq& \sum_{l|d} g(l)\sum_{h<\xi/d,\,dh|\Pi}
g(h)\\
&=&\sum_{l|d} g(l)\,G_{d}\big(\frac{\xi}{d},\,z\big)\;, 
\end{eqnarray*}
since $\xi/d\,\leq\,\xi/l$ and
$g(h)\,\geq\,0\,.$ The inequality (2.6) follows now from the
identity
$$\sum_{l|d} g(l)\,=\,\prod_{p|d}(1+g(p))\,=\,
\prod_{p|d}\big
\{1+\frac{\omega(p)}{p}\,\big(1-\frac{\omega(p)}{p}\big)^{-1}\big
\}=\prod_{p|d}\big(1-\frac{\omega(p)}{p}\big)^{-1}\;,$$ 
in view of the condition $(\Omega_{1})\;.$ \\
\bigskip

\noindent {\bf The ``Fundamental
Theorem for Selberg's sieve''.} \emph{Assume $(\Omega_{1})\,.$ Then
$$S(\cl{A},\cl{P};z)\,\leq\, \frac{X}{G(\xi,\,z)}\,+\,
\sum_{d<y,\,d|\Pi(\cl{P},\,z)} 3^{\nu(d)}\,|R_{d}|$$
for $\xi\geq 1\,,$ where $$\xi:=y^{\frac{1}{2}}\,, \quad
\nu(d):=\sum_{p|d,\,p\in \op}1\;,\quad
G(\xi,\,z):=\sum_{l|\Pi(\cl{P},z),\,l<\xi} g(l)$$ with
$$g(d)=\frac{\omega(d)}{d}\,\prod_{p|d}
\big(1-\frac{\omega(p)}{p}\big)^{-1}\,.$$}\\
{\bf Proof.} It follows from the relations (2.3) and (2.6) that
$|\lambda_{d}|\leq 1$. Therefore
$$ |\mu^{+}(d)|\,=\big|\sum_{[d_{1},\,d_{2}]=d} 
\lambda_{d_{1}}\lambda_{d_{2}}\big|\,\leq\,
\sum_{[d_{1},\,d_{2}]=d} 1 \,=\,\sum^{\nu(d)}_{\lambda=0} \twoar
{\nu(d)}{\lambda}\,2^{\lambda}\,=\,3^{\nu(d)}\,,$$  for square-free
$d$.
Consequently, the assertion of the theorem follows
from (2.5).

\renewcommand{\thesection}{3}
\section{Some applications}
We prove two corollaries of the  Fundamental Theorem.
\begin{cor}
Let $\cl{P}\,=\,\{p \in \op:\,p\nmid k\,\}$ and
let
$$\omega(p)=\twos{0,}{\quad p|k}{1,}{\quad p\nmid k}\;.$$
Then 
\begin{equation}
S(\cl{A},\cl{P};z)\leq \frac{k}{\varphi(k)}\,\frac{X}{\log z} +
\twosum{d|\Pi(\cl{P},z),\,d<z^{2}}{(d,\,k)=1}
3^{\nu(d)}\,|R_{d}|\;.
\end{equation}
\end{cor}
{\bf Proof.} When $l|\Pi(\cl{P},z)$ we have
$$g(l)=\frac{\omega(l)}{l}\,\prod_{p|l}\Bg{1-\frac{\omega(p)}{p}}^{-1}=
\frac{1}{\varphi(l)}\;$$ so that
$$G(z,\,z)=\sum_{l|\Pi(\cl{P},z),\,l<z}g(l)\,=
\twosum{(l,\,k)=1}{l<z}
\frac{\mu(l)^{2}}{\varphi(l)}\,=:\,H_{k}(z)\;.$$ Let  $K(n)$ be the
largest square-free ivisor of $n$.  One then obtains 
$$H_{1}(z)=\sum_{l<z} \frac{\mu(l)^{2}}{l}\,\prod_{p|l}
\Bg{1-\frac{1}{p}}^{-1}\,=
\sum_{l=p_{1}...p_{h}<z}\;\prod^{h}_{i=1}(p_{i}-1)^{-1}\,=$$
$$\twosum{\alpha_{i}\geq 1}{p_{1}...p_{h}<z}
\frac{1}{p^{\alpha_{1}}_{1}...p^{\alpha_{h}}_{h}}\,=\,
\sum_{K(n)<z} \frac{1}{n}\,\geq \,\sum_{n<z} \frac{1}{n}\,\geq\,\log z\,.$$
On the other hand, for any square-free $k$ we have
\begin{eqnarray*}
H_{1}(z)&=&\sum_{n<z}\frac{\mu(n)^{2}}{\varphi(n)}\\
&=&\sum_{l|k}\twosum{n<z}{(n,\,k)=l} \frac{\mu(n)^{2}}{\varphi(n)}\\
&=&\twosum{l|k,\;hl<z}{(h,\,k/l)=1}\frac{\mu(lh)^{2}}{\varphi(lh)}\\
&=&\sum_{l|k}\frac{\mu(l)^{2}}{\varphi(l)}\,\twosum{h<z/l}{(h,\,k)=1}
\frac{\mu(h)^{2}}{\varphi(h)\,}\\
&=&\sum_{l|k}\frac{\mu(l)^{2}}{\varphi(l)}\,H_{k}\Bg{\frac{z}{l}}\\
&\leq&H_{k}(z)\,\sum_{l|k}\frac{\mu(l)^{2}}{\varphi(l)}\\
&=&H_{k}(z)\prod_{p|k}\Bg{1+\frac{1}{p-1}}\\
&=&\frac{k}{\varphi(k)}\,H_{k}(z)\,.
\end{eqnarray*}
Combining the last two estimates, one obtains
\begin{equation}
H_{k}(z)\,\geq\,\frac{\varphi(k)}{k}\,H_{1}(z)\,\geq\,
\frac{\varphi(k)}{k}\,\log z\,.
\end{equation}
Thus 
\[G(z,z)\ge \frac{\phi(k)}{k}\log z\]
so that Corollary 3.1 follows from the Fundamental Theorem with  $\xi=z$.

\begin{cor} Let $\cl{P}\,=\,\{p\in \op:\,p\nmid k\,\}$ and let
$$\omega(p)=\twos{0,}{\quad p|k,}{\frac{p}{p-1},}{\quad p\nmid k.}$$ 
Then
$$S(\cl{A},\cl{P};z)\,\leq\,\Bg{\twoprod{p|k}{p\neq
2}\frac{p-1}{p-2}}\,C_{2}\,
\frac{X}{\log z}\,\left\{1+O\bg{\frac{1}{\log z}}\right\}
+\twosum{d<z^{2},\,d|\Pi(\cl{P},z)}{(d,\,k)=1}
3^{\nu(d)}\,|R_{d}|$$ for even $k$, where
$$C_{2}:=2\,\prod_{p\geq 3}{\Bg{1-\frac{1}{(p-1)^{2}}}}$$ is the ``twin 
prime constant''.
\end{cor}
{\bf Proof.} We have $g(p)=0 $ if $ p|k$, and
$$g(p)=\frac{\omega(p)}{p}\,\Bg{1-\frac{\omega(p)}{p}}^{-1}\,=\,
\frac{1}{(p-1)\,\Bg{1-\frac{1}{p-1}}}\,=\,\frac{1}{p-2}\,=\,
\frac{1}{\varphi(p)}\,(1+g(p))$$ 
otherwise. Therefore 
$$g(d)=\frac{1}{\varphi(d)}\sum_{l|d}g(l)\,\mu(l)^{2}$$
as soon as $(d,\,k)=1$ and $\mu(d)^{2}=1$.  Moreover, $g(d)=0$ when
$(d,\,k)\neq 1\;.$ Thus
\begin{eqnarray*}
G(z,\,z)&=&\sum_{d|\Pi(\cl{P},z),\,d<z}g(d)\\
&=&\sum_{d<z,\,(d,\,k)=1\;,l|d}
\frac{\mu(d)^{2}}{\varphi(d)}\,g(l)\,\mu(l)^{2}\\
&=&\twosum{(ml,\,k)=1}{ml<z}
\frac{\mu(ml)^{2}}{\varphi(ml)}\,g(l)\,\mu(l)^{2}\\
&=&\sum_{l<z,\,(l,\,k)=1} \frac{\mu(l)^{2}\,g(l)}{\varphi(l)}
\sum_{(m,\,lk)=1,\,m<\frac{z}{l}}\frac{\mu(m)^{2}}{\varphi(m)}\\
&=&\sum_{l<z,\,(l,\,k)=1}\frac{\mu(l)^{2}\,g(l)}{\varphi(l)}\,H_{lk}
\Bg{\frac{z}{l}}\,.
\end{eqnarray*}
In view of (3.2),
\begin{eqnarray*}
G(z,z)&\geq &\sum_{l<z,\,(l,k)=1}\frac{\mu(l)^{2}\,g(l)}{\varphi(l)}\,
\frac{\varphi(kl)}{kl}\,\log\bg{\frac{z}{l}}\\
&=&\frac{\varphi(k)}{k}\sum_{l<z,\,(l,k)=1}\frac{\mu(l)^{2}\,g(l)}
{l}\,\log\bg{\frac{z}{l}}\\
&\geq&\frac{\varphi(k)}{k}\twosum{l=1}{(l,k)=1}^{\infty}
\frac{\mu(l)^{2}\,g(l)}{l}\,
\log\bg{\frac{z}{l}}\\
&=&\frac{\varphi(k)}{k}(\log z)\,\prod_{p\nmid k}
\bg{1+\frac{1}{p\,(p-2)}}-\frac{\varphi(k)}{k}\twosum{l=1}{(l,k)=1}^{\infty}
\frac{\mu(l)^{2}\,g(l)}{l}\,\log l\\
&=&\frac{\varphi(k)}{k}\,\prod_{p\nmid k}\bg{1+\frac{1}{p\,(p-2)}}\,
\{\log z+O(1)\},
\end{eqnarray*}
so that
\[
\frac{1}{G(z,z)}\,\leq\,\frac{1}{\log z}\Bg{1+O\bg{\frac{1}{\log
z}}}\,\Bg{2\prod_{p|k}\frac{p-1}{p-2}}\, \prod_{p\geq
3}\bg{1-\frac{1}{(p-1)^{2}}}
\]
since $$\frac{\varphi(p)}{p}\,\bg{1+\frac{1}{p\,(p-2)}}^{-1}\,=\,
\frac{p-2}{p-1}$$
and
$$\bg{1+\frac{1}{\,p(p-2)}}^{-1}\,=\,\bg{1-\frac{1}{(p-1)^{2}}}$$
for $p\neq 2$.  Corollary 3.2 now follows from the Fundamental Theorem.

Our applications will require one further result.
\begin{lemma}
Let $h\in\N$ and set
$$
S_{1}=\sum_{d<x}{\mu(d)^{2}}\,h^{\nu(d)},\;\;\;\;
S_{2}=\sum_{d<x}\frac{\mu(d)^{2}}{d}\,h^{\nu(d)}\,.$$ We then have
\[
S_{1}\leq\,x\,(1+\log x)^{h}\quad \mbox{and} \quad
S_{2}\leq\,(1+\log x)^{h}\,.
\]
\end{lemma}
{\bf Proof.} Clearly,
$$S_{1}\leq\,\sum_{d<x}\mu(d)^{2}\,\frac{x}{d}\,h^{\nu(d)}\,=\,x\,S_{2}\,.$$
Moreover,
\begin{eqnarray*}
S_{2}&=& \sum_{d<x}\frac{\mu(d)^{2}}{d}\,\twosum{(d_{1},...,d_{h})}
{d_{1}...d_{h}=d} 1 \\
&\leq& \sum_{d\geq 1}\frac{1}{d}\twosum{d_{1}...d_{h}=d}{d_1,\ldots,d_h<x}
\mu(d_{1})^{2}...\mu(d_{h})^{2}\\
&=&\Bg{\sum_{d<x}\frac{\mu(d)^{2}}{d}}^{h}\leq\,(1+\log x)^{h},
\end{eqnarray*}
as asserted. Here we have used the fact that if $\mu(d)^{2}=1$ and 
$d=p_{1}...p_{h}$, then
$$h^{\nu}=\#\{(d_1,\ldots,d_h)\in\N^h:\,d_{1}...d_{h}=d\}.$$

{\bf First application.} Let
$\cl{A}=\{\,n\,|\,x<n\leq\,x+y\,\}$ and $\cl{P}=\op$, so that
$$\#\cl{A}_{d}=\frac{y}{d}\,+\,O(1),\;\;\;\omega(d)=1 \quad \mbox{for} 
\quad d \in \N\,.$$
Then
$$S(\cl{A},\cl{P};z)\,\leq\,\frac{y}{\log z}\,+\,O\bg{\sum_{d\leq
z^{2}} 
3^{\nu(d)}\,\mu(d)^{2}}$$
by Corollary 3.1.
In view of Lemma 3.1, this gives
$$S(\cl{A},\cl{P};z)\,\leq\,\frac{y}{\log z}\,+\,O\bg{z^{2}(1+\log z)^{3}}\,.$$
On the other hand,
$$S(\cl{A},\cl{P};z)\,=\#\{\,n\,:\,x<n\leq\,x+y,\; p|n\,\Rightarrow \,
p\geq z\,\}\,\geq\,\pi(x+y)\,-\,\pi(x)\,-\,z\,.$$
Taking $z=\frac{\sqrt{y}}{(\log y)^{3}}$, it follows that
\begin{equation}\pi(x+y)-\pi(x)\leq \frac{2y}{\log y}\,+\,
O\Bg{ \frac{y\,\log\log
y }{(\log y)^{2}} } \quad \mbox{for}\; x,y\geq 2\,.
\end{equation}
In particular,
$$\pi(y)\leq\frac{2y}{\log y}\,+\,O\Bg{\frac{y\,\log\log y}
{(\log y)^{2}}}\,.$$ It has been proved (Heath-Brown \cite{HB712}) that
\begin{equation}
\pi(x+y)-\pi(x)\,\sim \,\frac{y}{\log x}
\end{equation}
for $y\geq x^{\frac{7}{12}}$.  In contrast, (3.4) is false (Maier
\cite{Maier}) for
$y\,\asymp\,(\log x)^{A}$, for any constant $A$. Thus (3.3) is useful for
relatively small $y$. Montgomery and Vaughan \cite{MV} have removed
the error term above and proved
that
\[
\pi(x+y)-\pi(x)\leq\,\frac{2y}{\log y}
\quad \mbox{for}\;\;x\,,\,y\geq 2\,.
\]
It has been conjectured that $\pi(x+y)\leq\,\pi(x)+\pi(y)$ for all
$x,y\ge 2$.  However Hensley and Richards \cite{Hens} have proved that
this would be 
incompatible with the $k$-tuples prime conjecture.  It is not 
clear at the moment
whether the factor 2 in (3.3) may be replaced by a smaller 
number.  Indeed Erd\H{o}s apparently believed that the constant may be
taken as 1, while Selberg is reputed to have suggested that no
constant below 2 is admissible.

{\bf Second application.} Suppose that $(l,\,k)=1$ and let
$$\cl{A}=\{n\leq x:\;n\equiv l\mod{k}\,\},$$ 
$$\pi(x;k,l):=\#\{\,p\in\op:\,p\leq x,\,p\equiv l\mod{k}\,\},$$
$$\cl{P}=\{\,p\,\in\,\mathbb{P}:\,p\nmid k\,\}\,.$$ 
Clearly,
\begin{equation}
\pi(x;k,l)\leq\,S(\cl{A},\cl{P};z)+1+\frac{z}{k}\,.
\end{equation}
Moreover,
$$ \#\cl{A}_{d}=\frac{x}{k}\frac{\omega(d)}{d}\,+\,O(1)$$
with
$$\omega(d)=\twos{1,}{\;\mbox{if}\;\;(d,k)=1,}{0,}
{\;\mbox{if}\;\;(d,k)\neq 1,}$$ 
so that
\begin{eqnarray*}
S(\cl{A},\cl{P};z)&\leq&\frac{k}{\varphi(k)}\,\frac{x}{k}\,
\frac{1}{\log z}\,+\,\sum_{d\leq z^{2}}3^{\nu(d)}\,\mu(d)^{2}\\
&=&
\frac{x}{\varphi(k)\log z}\,+\,O\bg{z^{2}(\log z)^{3}}
\end{eqnarray*}
by Corollary 3.1 and Lemma 3.1. Let
$z=(x/k)^{1/2} (\log x/k)^{-3}$.
Assuming that $x\ge 4k$, say, the estimate (3.5) yields the
following result.

\noindent {\bf The Brun-Titchmarsh Theorem.}  {\em We have}
\begin{equation}
\pi(x;k,l)\leq\,\frac{2x}{\varphi(k)\log x/k}\,
+\,O\Bg{\frac{x}{\varphi(k)}\,
\frac{\log\log x/k}{(\log x/k)^{2}}}\,.
\end{equation}
{\em for} $x\ge 4\varphi(k)$.

The Siegel-Walfisz Theorem gives 
$$\pi(x;k,l)\sim \frac{x}{\varphi(k)\log x} \quad \mbox{for} \quad k
\ll(\log x)^{A}\,$$
for any fixed $A$, so that the constant 2 in the Brun-Titchmarsh
Theorem may be replaced by 1 if $x$ is sufficiently large compared
with $x$.  Moreover Montgomery and Vaughan \cite{MV} have proved that
\[
\pi(x;k,l)\leq \frac{2x}{\varphi(k)\log x/k}
\quad \mbox{for}  \quad k<x.
\]
The constant 2 in (3.6) is presumably hard to improve, for it is known
that if one could replace $2$ by $2-\delta$ with a positive constant 
$\delta$, then
it would follow that there are no ``Siegel-Landau zeros". 

{\bf Third application.} Let
$$\cl{A}=\{\,2N-p\,:\,p\in\,\cl{P},\;3\leq p\leq 2N-3\,\},\;\;\;
\cl{P}=\{\,p\in\,\op:\;(p, 2N) = 1\,\}\;,$$ 
so that
\begin{eqnarray*}\cl{A}_{d}&=&
\#\{\,p\in\,\op:\;p\equiv 2N\mod{d},\,p\nmid 2N,
\;p\leq2N-3\,\}\\
&=&\pi(2N;d,2N)+O\big(1+\nu(N)\big)\,.
\end{eqnarray*}
When $\text{h.c.f.}(d,2N)=1$ we expect that
$$\pi(N;d,2N)\sim \frac{\li(2N)}{\varphi(d)}$$
for $N$ large compared with $d$.  We therefore take $X=\li(2N)$ and 
$\omega(d)=\varphi(d)^{-1}d$. 
On writing
$$E(x;k,l)=\pi(x;k,l)-\frac{\li(x)}{\varphi(k)}\,$$
it follows that
$$R_{d}=E(2N;d,2N)\,+\,O\big(1+\nu(N)\big)\,.$$ 
On the other hand
$$r(2N):=\#\{\,(p,q)\,|\,p,q\in\,\op,\;p+q=2N\,\}\leq 
S(\cl{A},\cl{P};z)\,+\,z\,$$
and we deduce from Corollary 3.2 that
$$r(2N)\leq \Bg{\twoprod{p|2N}{p\neq2}\frac{p-1}{p-2}}\,C_{2}\,
\frac{\li(2N)}{\log
z}\,\Bg{1+O\bg{\frac{1}{\log z}}}\,+\,z$$
\[\hspace{6mm}
+\,O\big(\sum_{d\leq
z^{2}}3^{\nu(d)}\,\mu(d)^{2}\,\{|E(2N;d,2N)|\,+\,1+\nu(N)\}\big)\,.
\]
In order to estimate the remainder sum 
\[S:=\sum_{d\leq z^{2}} 3^{\nu(d)}\mu(d)^{2}E_{1}(2N,d),\]
where
$$E_{1}(x,k):=\max_{y\leq x,\,(l,k)=1}|E(y;k,l)|\,,$$
we shall use the following
well-known result.

\noindent {\bf The Bombieri-Vinogradov Theorem.} {\it For every $c_{1}>0$, 
there is a positive constant $c_{2}$ such that} 
$$\sum_{k\leq x^{\frac{1}{2}}\,(\log x)^{-c_{2}}}E_{1}(x,k)\ll
\frac{x}{(\log x)^{c_{1}}}.$$

Using Cauchy's inequality with Lemma 3.1 we find that
\begin{eqnarray*}
S^{2}&\leq &\big\{\sum_{d<z^{2}}(3^{\nu(d)}\mu(d)^{2})^{2}\,
\frac{1}{d}\big\}\,\big\{\sum_{d\leq z^{2}}d\,E_{1}(2N,d)^{2}\big\}\\
&\ll & (\log z)^{9}\sum_{d<z^{2}}d\,E_{1}(2N,d)^{2}.
\end{eqnarray*}
It is trivial that
\[\pi(x;k,l)\ll\frac{x}{k}\ll\frac{x}{\varphi(k)}\]
for $k\le x$, so that
\[E(x;k,l)\ll\pi(x;k,l)+\frac{x}{\varphi(k)}\ll\frac{x}{\varphi(k)}\ll
\frac{x\log x}{k}.\]
Thus
\[\sum_{d<z^{2}}d\,E_{1}(2N,d)^{2}\ll N(\log N)\sum_{d<z^{2}}E_{1}(2N,d)
\ll N^2(\log N)^{-15},\]
by the Bombieri-Vinogradov Theorem with $c_1=16$, on taking 
$$z^{2}=\sqrt{N}\,(\log N)^{-c_{2}}.$$
We therefore conclude that
$$S^{2}\ll\frac{N^{2}}
{(\log N)^{6}}.$$
Thus
\[
\sum_{d\leq z^{2}}3^{\nu(d)}\mu(d)^{2}E_{1}(2N,d)=O(N\,(\log
N)^{-3})\,.
\]
Moreover, it follows from Lemma 3.1 that
\[
\big(1+\nu(N)\big)\sum_{d\leq z^{2}}3^{\nu(d)}\mu(d)^{2}\ll
\big(1+\nu(N)\big)\sqrt{N}(\log N)^{3}\ll \frac{N}{(\log N)^{3}}\,.
\]
We may now deduce the following result from (3.1).
\begin{theorem}
We have 
$$r(2N)\leq \{4+O(\frac{\log\log N}{\log N})\}a(N),$$ 
with
$$a(N)=\Bg{\twoprod{p|2N}{p\neq 2}\frac{p-1}{p-2}}\,C_{2}\,
\frac{2N}{(\log N)^{2}}.$$
\end{theorem}
It is conjectured that one may improve the Bombieri-Vinogradov Theorem
to say that for any $\ep>0$ and any $c_1>0$ one has
$$\sum_{k\leq x^{1-\ep}}E_{1}(x,k)\ll
\frac{x}{(\log x)^{c_{1}}}\,.$$
One would then obtain a bound
\[r(2N)\leq (2+o(1))a(N)\]
in a completely analogous fashion.  However the best unconditional
result is due to Chen \cite{Chen7}, in which the constant 4 is
reduced to $3.9171$. For comparison we note that 
it is conjectured that $r(2N)\sim a(N)$.

The following theorem can be proved in the same way as Theorem 3.1 (Exercise!).
\begin{theorem}
For any positive integer $k$ we have
$$\#\{\,p\le x\,:\,p,p+2k\in \op\}\leq
4\bg{\prod_{p|2k} \frac{p-1}{p-2}}\,C_{2}\,\frac{x}{(\log x)^{2}}\, 
\bg{1+O\bg{\frac{\log\log x}{\log x}}}\,.$$
\end{theorem}
\begin{cor} ({\it Viggo Brun})
\par We have 
$$\sum_{p,\,p+2\in\op}\frac{1}{p}\,<\,\infty.$$
\end{cor}
Bombieri, Friedlander and Iwaniec \cite{BFI} have proved a variant of 
Theorem 3.2 
with the constant 4 replaced by $7/2$. Their method does not establish
a result uniform in $k$ and is therefore not applicable to Theorem
3.1.  More complicated methods allow
one to
reduce the constant in Theorem 3.2 further slightly.
\bigskip

We proceed to discuss briefly some other applications of Selberg' sieve.
{\bf Definition.} {\it Suppose that
\begin{equation}
\sum_{w\leq p<z}\frac{\omega(p)\,\log p}{p}\, =\,\kappa\,
\log(\frac{z}{w})\,+\,O(1)\quad \mbox{for} \quad 2\leq w\leq
z.
\end{equation}
Then the constant $\kappa$ is called the {\bf dimension} of the sieve problem.}

One can get by with a slightly weaker assumption in fact.  The above 
definition corresponds to a version of the condition $\Omega(\kappa,L)$
in the book
by Halberstam and Richert \cite[page 142]{HR}.

{\bf Remark.} Since
$$\sum_{w\leq p<z}\frac{\log p}{p}\,=\,\log(\frac{z}{w})\,+\,O(1)\,,$$
the dimension of the sieve problem coincides with the ``average
value" of $\omega(p)$. Note that in the two cases considered in 
this section we have $\kappa = 1$.  This is clear for both
$$\omega(p)\,=\,\twos{1,}{\quad p\nmid k,}{0,}{\quad p|k,}$$ 
and
$$\omega(p)\,=\,\twos{1+\frac{1}{p-1},}{\quad p\nmid k,}{0,}{\quad p|k,}$$

In general, for the sieve problem of dimension $\kappa$, one obtains
\begin{eqnarray*}
G(z,z)&=&\twosum{d|\Pi(\cl{P},z)}{d\leq z}\frac{\omega(d)}{d}\,
\prod_{p|d}\bg{1-\frac{\omega(p)}{p}}^{-1}\\
&=&\sum_{d\leq z} \frac{\omega(d)\mu(d)^{2}}{d}\,\prod_{p|d}
\bg{1-\frac{\omega(p)}{p}}^{-1}\,.
\end{eqnarray*}
Though we shall not prove it, it turns out that
$$G(z,z)=\frac{1}
{e^{\gamma \kappa}\,\Gamma(\kappa+1)}\,\prod_{p<z}
\bg{1-\frac{\omega(p)}{p}}^{-1}\,\bg{1+O\bg{\frac{1}{\log z}}}$$
see Halberstam and Richert \cite[(5.3.1)]{HR}. It
therefore  follows from the Fundamental Theorem that
\begin{eqnarray}
S(\cl{A},\cl{P};z)&\leq &X\,\prod_{p<z}\bg{1-\frac{\omega(p)}{p}}\, 
e^{\gamma \kappa}\,\Gamma(\kappa+1)\,\bg{1+O\bg{\frac{1}{\log
z}}}\nonumber\\
&&\hspace{1cm}+\sum_{d<z^{2},\,d|\Pi(\cl{P},z)}3^{\nu(d)}|R_{d}|.
\end{eqnarray}
Note that if $\omega(p)=1$, then
$$\prod_{p<z}\bg{1-\frac{1}{p}}\sim\frac{e^{-\gamma}}{\log z}\,.$$ 
One therefore obtains
$$S(\cl{A},\op;z)\leq\frac{X}{\log z}\,\bg{1+O\bg{\frac{1}
{\log z}}}\,+\,\sum_{d<z^{2},\,d|\Pi(\op,z)}3^{\nu(d)}|R_{d}|\,,$$
cf. (3.1).

{\bf Example.} Let 
$$\cl{A}=\{\,n(2N-n)\,:\,3\leq n\leq 2N-3\,\},$$
so that
\begin{eqnarray*}
\#\cl{A}_{d}&=&\#\{\,n\,|\,1\leq n\leq d,\;d|\,n(2N-n)\,\}\,
\bg{\frac{2N}{d}+O(1)}\\
&=&\frac{\omega(d)}{d}\,2N +
O(\omega(d))\,,
\end{eqnarray*}
where 
$$\omega(p)=\twos{1,}{\quad p|\,2N,}{2,}{\quad p\nmid 2N.}$$ 
In this case $\kappa=2$ and the estimate (3.8) implies $r(2N)\leq 
(8+o(1))a(N)$
( Exercise!).  This should be compared with Theorem 3.1, which
was deduced using a sieve of dimension 1, together with the 
Bombieri-Vinogradov Theorem.

\renewcommand{\thesection}{4}
\section{The parity phenomenon and limitations to sieve
methods}

The optimisation problem for the upper bound sieve amounts to the question of 
minimising the linear functional
$$l(\mu^{*}):=X\sum_{d|\Pi(\cl{P},z)}\frac{\mu^{*}(d)\,\omega(d)}{d}\,+\,
\sum_{d|\Pi(\cl{P},z)}\mu^{*}(d)R_{d}$$
under the additional condition that
$$\sum_{d|n,\,d|\Pi(\cl{P},z)}\mu^{*}(d)\geq \twos{1,}{
(n,\Pi(\cl{P},z))=1,}{0,}{\mbox{otherwise,}}$$
(cf. (1.5)--(1.7)).  We can view this as a linear programming
problem. In the standard formulation of a linear programing problem
one takes a real $m\times n$ matrix $A$ and real column vectors
$\b{b}$ and $\b{c}$ of lengths $m$ and $n$ respectively.  One then
seeks to minimize $\b{c}^t\b{x}$ over all column vectors
$\b{x}\in\R^n$, subject to the conditions that $\b{x}\ge\b{0}$ and
$A\b{x}\ge\b{b}$.  (Here $\b{z}\ge\b{w}$ means that $z_i\ge w_i$ for
each index $i$.)  In our problem the vector of values of $\mu^*(d)$ is
not required to be non-negative, so we introduce two new functions
$\mu^*_{\pm}(d)$ with $\mu^*_{\pm}(d)\ge 0$ for all $d$, and write
$\mu^*(d)=\mu^*_+(d)-\mu^*_-(d)$.  We can then produce a linear
programming problem in standard form.

A great deal can be learnt about a linear programing problem by
studying its ``dual''.  For the problem described above, the dual
problem is to maximize $\b{y}^t\b{b}$ over column vectors
$\b{y}\in\R^m$, subject to $\b{y}\ge\b{0}$ and $\b{y}^tA\le\b{c}^t$.
Under these constraints one clearly has
\begin{equation}
\b{c}^{t}\b{x}\geq \b{y}^{t}A\b{x}\geq \b{y}^{t}\b{b}\,
\end{equation}
and the Duality Theorem states that there exist vectors $\b{x},\b{y}$
for which equality is attained.

In the context of the upper bound sieve problem, it transpires that
the dual problem is essentially that of finding a sequence $\cl{A}$,
with prescribed function $\omega(d)$, and with suitably small
remainders $R_d$, for which $S(\cl{A},\cl{P};z)$ is as large as possible.
We can interpret the inequalities (4.1) as saying that, for any vector
$\b{y}=\b{y}_0$ which satisfies the relevant constraints, we must have
\[\inf \b{c}^{t}\b{x}\geq \b{y}_0^{t}\b{b},\]
and moreover, if we have vectors $\b{x}$ and $\b{y}$ which both satisfy
the relevant constraints, and for which
$\b{c}^{t}\b{x}=\b{y}^{t}\b{b}$, then both $\b{x}$ and $\b{y}$ must be
extremal.  

For the upper bound sieve problem, any sequence $\cl{A}$ defining a
problem of dimension 1 will therefore produce a lower bound on the
possible values of $l(\mu^*)$.  Moreover if we can find a sequence
$\cl{A}$ and
a set of coefficients $\mu^*(d)$ for which $S(\cl{A},\cl{P};z)$
(corresponding to $\b{y}^{t}\b{b}$) and $l(\mu^*)$ (corresponding to 
$\b{c}^{t}\b{x}$) are approximately equal, then both must be
essentially optimal.  

We therefore examine in detail the following two
sequences, first discussed by Selberg.  Let $\Omega(n)$ be the number
of prime factors of $n$, counted according to multiplicity, and define
the {\em Liouville function} $\lambda(n)$ by
\[\lambda(n)=(-1)^{\Omega(n)}.\]
We then set
$$\cl{A}^{+}=\{n\in\,\N:\,n\leq x,\,\lambda(n)=-1\,\}\,,$$
and
$$\cl{A}^{-}=\{n\in\,\N:\,n\leq x,\,\lambda(n)=+1\,\}\,,$$
which will relate to the upper bound and lower bound problems
respectively.  In the case of the first sequence we have
$$S(\cl{A}^{+},\op;z)=\pi(x)-\pi(z)=\frac{x}{\log
x}\,+\,O\bg{\frac{x}{(\log x)^{2}}}\,+\,O(z)
\quad \mbox{for}\; z> x^{\frac{1}{3}}\,.$$ 
Let us see how this compares with the bound given by the Selberg sieve.
To bound $R_d$ for the sequences $\cl{A}^{\pm}$ we note that if
\[L(x):=\sum_{n\le x}\lambda(n),\]
then $L(x)\ll E_c(x)$, where $c$ is a suitable positive constant and 
\[E_c(x):=\exp(-c\sqrt{\log x}).\]
(This follows by a similar analysis to that used for the summatory
function of $\mu(n)$.)  Now, if we let $X=\frac{x}{2}$ then
\begin{eqnarray*}
\#\cl{A}^{\pm}_{d}&=&\#\{m\leq \frac{x}{d}:\;\lambda(md)=
\mp1\,\}\\
&=&\#\{m\leq
\frac{x}{d}:\;\lambda(m)=\mp\lambda(d)\,\}\\
&=&\frac{[x/d]}{2}\mp\lambda(d)\frac{L(x/d)}{2}\\
&=&\frac{X}{d}\,+\,O\bg{E_{c}\bg{\frac{x}{d}}}\,,
\end{eqnarray*} 
for $d\le x$, and hence the remainder sum in the Fundamental Theorem is
\begin{eqnarray*}
\sum_{d<y,\,d|\Pi(\cl{P},z)}3^{\nu(d)}\,|R_{d}|&\ll&
\sum_{d<y,}3^{\nu(d)}\,E_{c}\bg{\frac{x}{d}}\\
&\leq & \Bg{\sum_{d<y,}9^{\nu(d)}d^{-1}}^{\frac{1}{2}}\, 
\Bg{\sum_{d<y}d\,E_{c}\bg{\frac{x}{d}}^{2}}^{\frac{1}{2}}\\
&\ll&
(\log y)^{\frac{9}{2}}\,\Bg{\sum_{d<y}d\,E_{c} 
\bg{\frac{x}{d}}^{2}}^{\frac{1}{2}}\,
\end{eqnarray*}
in view of Lemma 3.3.  Furthermore,
\begin{eqnarray*}
\sum_{d<y}dE_{c}\bg{\frac{x}{d}}^2&=&x^2\sum_{d<y}d^{-1}
\exp(-2c\sqrt{\log x/d})\\
&\le& x^2\sum_{d<y}d^{-1}\exp(-2c\sqrt{\log x/y})\\
&\le& x^2(\log y)\exp(-2c\sqrt{\log x/y})
\end{eqnarray*}
for $y<x\,.$ Thus
\begin{equation}
\sum_{d<y}3^{\nu(d)}\,|R_{d}|\ll
(\log x)^{5}\,x\,e^{-c\,\sqrt{\log x/y}}\ll
x\,(\log x)^{-2}
\end{equation}
for $y\le E_{1}(x)$, say. By Corollary 3.1, we have
$$S(\cl{A}^{+},\op;z)\leq \frac{X}{\log z}\,+\,
\sum_{d<z^{2}} 3^{\nu(d)}|R_{d}|.$$
On taking
$y=E_1(x)$ and $z=y^{1/2}$ one deduces from 
these estimates on recalling that $x=2X$, that
\[
S(\cl{A}^{+},\op;z)\leq \frac{X}{\log(\sqrt{E_{1}(x)})}\, +\,
O\bg{\frac{x}{(\log x)^{2}}}\,=\frac{x}{\log x}\,+ \,
O\bg{\frac{x}{(\log x)^{3/2}}}\,.
\]
Since
$$S(\cl{A}^{+},\op;z)=\pi(x)-\pi(z)=\frac{x}{\log
x}\,+\,O\bg{\frac{x}{(\log x)^{2}}}\,+\,O(z)
\quad \mbox{for}\; z> x^{1/3}\,,$$ 
we conclude that the estimate (3.1) cannot be improved on for $z>
x^{1/3}$. 

Thus the Selberg sieve is best possible in this situation.  Our
remarks about linear programming then show that the Selberg sieve
coefficients are an essentially optimal solution to the minimiztion
problem for $l(\mu^*)$, and that the sequence $\cl{A}^+$ is a
corresponding solution for the dual problem.

Turning to the lower bound sieve problem, we see that we can satisfy
the relevant constraints
$$\sum_{d|n,\,d|\Pi(\cl{P},z)}\mu^{*}(d)\leq \twos{1,}{\quad \mbox{if}
\quad
(n,\Pi(\cl{P},z))=1,}{0,}{\quad\mbox{otherwise,}}$$
by taking $\mu^*(d)$ to be identically zero.  For this choice we
produce the trivial lower bound
\[S(\cl{A}, \cl{P}; z) \ge X \sum_{d|\Pi (\cl{P}, z)} \frac{\mu^*(d)
\omega(d)}{d} - \sum_{d|\Pi (\cl{P}, z)} |\mu^*(d)|\, |R_{d}|=0.\]
We now observe that for the sequence $\cl{A}^-$ we have
$S(\cl{A}^-,\op,z)=1$ for $z>x^{1/2}$, since only the integer
$1\in\cl{A}^-$ is counted.  Thus the coefficients $\mu^*(d)=0$ are
essentially best possible for the linear programming problem in this
situation, and the sequence $\cl{A}^-$ is essentially optimal for the
corresponding dual problem.

Thus no set of lower bound sieve coefficients $\mu^*(d)$ with
$|\mu^*(d)|\le 3^{\nu(d)}$ can produce
\[\sum_{d|\Pi(\op, z),\,d<y} \frac{\mu^*(d)}{d}\gg\frac{1}{\log y}\]
when $z> x^{1/2}$.  In particular one cannot show that $\pi(x)\gg
x/\log x$ by sieve methods alone.

The two sequences $\cl{A}^+$ and $\cl{A}^-$ produce the same
information input for the sieve.  They have the same $X$, the same
function $\omega(d)$, and their remainders $R_d$ have the same order
of magnitude.  Thus there is no way that the sieve machinery can
distinguish them.  It is for this reason that the sieve encounters the
{\em parity phenomenon}, since it is unable to distinguish integers
for which $\Omega(n)$ is even, from those for which $\omega(n)$ is odd.

The sequences $\cl{A}^{\pm}$ have been shown to be essentially
extremal for $z> x^{1/2}$, but it transpires that they are optimal for
all $z$.  To examine this fact we define
\begin{eqnarray*}
S^{\pm}(x,s):&=&S\bg{\cl{A}^{\pm},\op;x^{1/s }}\\
&=&\#\{n\in\N:\,n\leq x,\,\lambda(n)=\mp1,\,
p|n\,\Rightarrow\,p\geq x^{1/s }\,\}\,
\end{eqnarray*}
for $s\ge 1$.  We classify the integers $n$ according to their
smallest prime factor $p$.  Then if
\[\eta_+=0,\;\;\;\eta_-=1\]
it follows that
\begin{eqnarray*}
S^{\pm}(x,s)&=&\sum_{p\geq x^{1/s }}\#\{m:\;
pm\leq x,\,\lambda(m)=\pm 1,\,p'|m\,\Rightarrow \,p'\geq p\,\}
\,+\,\eta_{\pm}\\
&=&\sum_{x^{1/s }\leq p\leq x} S^{\mp}\bg{\frac{x}{p},\,
\frac{\log x}{\log p}-1}\,+\,\eta_{\pm},
\end{eqnarray*}
since 
\[p=\left(\frac{x}{p}\right)^{1/s'}\]
with
\[s'=\frac{\log x/p}{\log p}.\]
This leads to recursion formulae for $S^{\pm}(x,s)$.  To produce
appropriate starting values for the recursions we note that
\[S^{-}(x,s)=O(1),\;\;\;(1\leq s\leq 2)\]
and 
\[S^{+}(x,s)=\pi(x)-
\pi(x^{1/s })+O(1)\;\;\;(1\leq s\leq 3).\]

Let us define continuous functions 
$F,f\,:\,\R_{>0}\,\rightarrow\,\R$ by
the relations 
$$F(s)=\frac{2\,e^{\gamma}}{s}\quad \mbox{for}\; 0<s\leq 3,\quad 
\quad f(s)=0\quad \mbox{for}\; 0<s\leq 2,$$
and
$$(s\,F(s))'=f(s-1),\quad (s\,f(s))'=F(s-1)\quad \mbox{for}\;s>2\,.$$
We note that these definitions show that
\[f(s)=2e^{\gamma}\frac{\log(s-1)}{s}\quad \mbox{for}\;\;2\leq s\leq 4.\]
We shall now prove the following estimates.
\begin{theorem}
Let $N\in\,\mathbb{N}$. Then we have
\begin{equation}
S^{+}(x,s)=\frac{x/2}{e^{\gamma}\,\log(x^{1/s })}\,
F(s)\,+\,O_{N}\bg{\frac{x}{(s-1)(\log x)^{2}}}
\end{equation}
and
\begin{equation}
S^{-}(x,s)=\frac{x/2}{e^{\gamma}\,\log(x^{1/s })}\,
f(s)\,+\,O_{N}\bg{\frac{x}{(\log x)^{2}}}
\end{equation}
for $1\leq s\leq N$.
\end{theorem}
{\bf Proof.} Since
$$\frac{x/2}{e^{\gamma}\,\log(x^{1/s })}\,F(s)=
\frac{x}{\log x}\;\;\;\;(0<s\le 3)$$
and 
$$S^{+}(x,s)=\pi(x)-\pi(x^{1/s })+O(1)\;\;\;\;(0<s\le 3)\,,$$
relation (4.3) holds for $N=2$ and $N=3$.  Similarly, when $N=2$ the
equation (4.4) follows from the facts that $f(s)=0$ and $S^-(x,s)=O(1)$
whenever $1\le s\le 2$.  We now prove (4.3) and (4.4) by induction on
$N$.  We shall consider only $S^{+}(x,s)$, leaving the discussion of 
$S^{-}(x,s)$ as an exercise.  Thus we assume that (4.3) and (4.4) hold
for $N$, and deduce that (4.3) holds for $N+1$.  We therefore let 
$N<s\leq N+1,$ with $N\geq 3$. Since
$$0<\frac{\log x}{\log p}-1\leq 2\quad  \mbox{for}\quad p\geq x^{1/s}\,,$$ 
it follows that
\begin{eqnarray*}
S^{+}(x,s)&=&\sum_{x^{1/s }\leq p\leq x}S^{-}\bg{\frac{x}{p},\,
\frac{\log x}{\log p}-1}\\
&=&\sum_{x^{1/3}< p\le x}S^{-}\bg{\frac{x}{p},\,\frac{\log x}{p}-1}
+\sum_{x^{1/s }\leq
p\le x^{1/3 }}S^{-}\bg{\frac{x}{p},\,\frac{\log x}{\log p}-1}\\
&=&\pi(x)-\pi(x^{1/3 })+O(1)+\sum_{x^{1/s }\leq
p<x^{1/3 }}S^{-}\bg{\frac{x}{p},\,\frac{\log x}{\log p}-1}.
\end{eqnarray*}
Moreover, if $p\geq x^{1/s }$ then $\frac{\log x}{\log p}-1
\leq s-1 \leq N $. Therefore, by the inductive assumption, we have
\begin{equation}
S^{-}\bg{\frac{x}{p},\,\frac{\log x}{\log p}-1}=
\frac{(2e^{\gamma})^{-1}x}
{p\,\log p}\, f\bg{\frac{\log x }{\log p}-1}\,
+\,O_{N}\bg{\frac{x}{p\,(\log x)^{2}}}
\end{equation}
for $x^{1/s }\le p\le x^{1/3}$.  By partial summation we find that
\begin{eqnarray*}
\sum_{x^{1/s }\le p\le x^{1/3}}\frac{1}
{p\,\log p}\, f\bg{\frac{\log x }{\log p}-1}&=&
\int^{x^{1/3 }}_{x^{1/s }}\frac{1}{t\log t}
f\bg{\frac{\log x}{\log t}-1}\frac{dt}{\log t}\\
&&\hspace{2cm}
+O_{N}\bg{\frac{x}{(\log x)^{2}}}\\
&=&\frac{1}{\log x}\int^{s}_{3}f(v-1)\,dv
+O_{N}\bg{\frac{x}{(\log x)^{2}}}\\
&=&\frac{1}{\log x}(sF(s)-3F(3))+O_{N}\bg{\frac{x}{(\log x)^{2}}},
\end{eqnarray*}
on substituting $t=x^{1/v}$.  To handle the error term in (4.5) we
note that
\[\sum_{x^{1/s }\le p\le x^{1/3}}\frac{x}{p(\log
  x)^2}\ll_N\frac{x}{(\log x)^2},\]
whence we conclude that
\begin{eqnarray*}
S^{+}(x,s)&=&\frac{x}{2\,e^{\gamma}\log x}\,
(sF(s)-3F(3))+\frac{x}{\log x}+O_{N}\bg{\frac{x}{(\log x)^{2}}}\\
&=&\frac{x}{2\,e^{\gamma}\log x^{1/s}}F(s)+O_{N}\bg{\frac{x}{(\log x)^{2}}},
\end{eqnarray*}
as required. 

{\bf Remarks.} \par 1) The properties
of the functions $f,\,F$ and their generalisations are discussed in detail
in the books by  Greaves \cite{GreavesBook} and Halberstam and Richert
\cite[Chapter 8]{HR}. 
\par 2) If we set
\[W(z):=\prod_{p\leq z}\bg{1-\frac{1}{p}},\]
then the Mertens formula (1.4) gives
$$\frac{x/2}{\,e^{\gamma}\log(x^{1/s })}=XW(x^{1/s
})+O_N(\frac{x}{(\log x)^{2}}),$$
so that (4.3) and (4.4) imply that
\[S^{+}(x,s)=X\,W(x^{1/s })F(s)\,+\,O_{N}\bg{\frac{x}{(s-1)(\log x)^{2}}}\]
and
\[S^{-}(x,s)=X\,W(x^{1/s })f(s)\,+\,O_{N}\bg{\frac{x}{(\log x)^{2}}}\]
respectively.

\renewcommand{\thesection}{5}
\section{The Rosser sieve}

A {\it combinatorial sieve} is defined by  choosing sets 
\[T^+(y),T^-(y)\subseteq\{d\in\N: d<y,\,\mu(d)^2=1\}\]
and taking
\begin{equation}
\mu^{\pm}(d)=\twos{\mu(d),}{\quad \mbox{if}
\;d\in\,T^{\pm}(y),}{0,}{\quad \mbox{otherwise.}}
\end{equation}
The sets $T^{\pm}(y)$ have to be chosen so that $\mu^{\pm}(d)$ satisfy
(1.5),  and (1.6) respectively.
As in the proof of (1.7) we have
\begin{equation}
S(\cl{A},\cl{P};z)\leq X\,\sum_{d|\Pi(\cl{P},z)}
\frac{\mu^{+}(d)\omega(d)}{d}\,
+\,\sum_{d|\Pi(\cl{P},z)}\mu^{+}(d)\,R_{d}\,.
\end{equation}
It follows from (5.1) and (5.2) that
$$|\sum_{d|\Pi(\cl{P},z)}\mu^{+}(d)\,R_{d}|\leq 
\sum_{d\leq y,\,d|\Pi(\cl{P},z)}|R_{d}|\,.$$ 
We have to choose $T^{+}(y)$ so as to optimise the main term in 
(5.2) subject to the condition (1.6). It follows from (1.8) and its 
analogue for $\mu^{-}$ that
\[S(\cl{A},\cl{P};z)\geq X\,\twosum{d|\Pi(\cl{P},z)}{d\in\;T^{-}(y)}
\frac{\mu^{-}(d)\omega(d)}{d}-\sum_{d\leq
  y,\,d|\Pi(\cl{P},z)}|R_{d}|\]
and
\[ S(\cl{A},\cl{P};z)\leq
 X\,\twosum{d|\Pi(\cl{P},z)}{d\in\;T^{+}(y)}\frac{\mu^{+}(d)\omega(d)}{d}
+\sum_{d\leq y,\,d|\Pi(\cl{P},z)}|R_{d}|\,.\]
One way to arrange for (1.5) and (1.6) to hold is as follows.  Let
$$T:=\{\,d\,|\,d\in\,\N,\;\mu(d)^{2}=1\,\},$$ 
write
$d=p_{1}p_{2}p_{3}...$ with $p_{1}>p_{2}>p_{3}...$ for $d\in\,T$, and let
$$T_{r}=\{\,d\,|\,d\in\,T,\;(\nu(d)<r)\;\mbox{or}\;(\nu(d)\geq r\,\&\,
{\mathfrak P}(p_{1},...,p_{r}))\,\}$$ for some predicate 
${\mathfrak P}$ to be defined later. Let
$$T^{+}(y)=\bigcap^{\infty}_{u=1}\,T_{2u-1},\quad T^{-}(y)=
\bigcap^{\infty}_{u=1}\,T_{2u}\,.$$
With these definitions we have the following lemma.
\begin{lemma} For $m\in\,\N$ and $T^{\pm}(y)$ defined as above, we have
$$\sum_{d|m}\mu^{-}(d)\leq \sum_{d|m}\mu(d)\leq \sum_{d|m}\mu^{+}(d)\,.$$
\end{lemma}
{\bf Proof.} Let
$$B_{2u-1}=(T\,\setminus\,T_{2u-1})\,\cap\,(\bigcap_{v<u}T_{2v-1})\,,$$
then 
$$T\,\setminus\,T^{+}(y)=\{d\in T:\,
\exists v\,(d\not\in \,T_{2v-1})\,\} =\bigcup^{\infty}_{u=1}\,B_{2u-1}\,.$$
Moreover we have $B_{j}\cap B_{k}=\emptyset$ if $j\neq k$. Therefore
$$\sum_{d|m}\mu(d)=\sum_{d|m,\,d\in\,T^{+}(y)}\mu(d)+U,$$
where
$$U:=\sum_{u=1}^{\infty}\;\sum_{d|m,\,d\in\,B_{2u-1}}\mu(d)\,.$$ 
We set
$$C_{r}=\{d\in\,B_{r}:\;\nu(d)=r\,\},$$
and if $d>1$ we write $p(d)$ for the smallest prime factor of $d$.  We
then define
$$Q(d)=\prod_{p\in\,\op,\,p<p(d)}p,\;\;\;(d>1).$$ 
Now, if $d\in\,B_{r}$, then $d\,\not\in\,T_{r}$. Thus $\nu(d)\geq r$ and
${\mathfrak P}(p_{1},...,p_{r})$ does not hold.
Write $d=ef$ with $e=p_{1}...p_{r}$ and $f|Q(e)$. Since $d\in\,B_{r}$,
it  follows
that each property ${\mathfrak P}(p_{1},...,p_{r-2})$, ${\mathfrak
  P}(p_{1},...,p_{r-4})$, $\ldots$ holds.  Now let
$$C_{r}=\{d\in B_{r}:\;\nu(d)=r\,\}\,,$$
so that $e\in\,C_{r}$.
Hence if $d\,\in\,B_{r}$ then we can write $d=ef$ with $e\in\,C_{r}$
and $f|Q(e)$.
Clearly, the decomposition $d=ef$ with $e\in\,C_{r}$ and $f|Q(e)$ is
unique. Conversely, if $d=ef$ with $e\in\,C_{r}$ and $f|Q(e)$, then
$d\in\,B_{r}$.
Therefore
\begin{eqnarray*}
\sum_{d|m,\,d\in\,B_{2u-1}}\mu(d)&=&
\twosum{ef|m,\,e\in C_{2u-1}}{f|Q(e)}\mu(ef)\\
&=&\twosum{e\in\,C_{2u-1}}{e|m}\mu(e)\,\sum_{f|m,\,f|Q(e)}\mu(f)\\
&=&\twosum{e|m,\,(m,Q(e))=1}{e\in\,C_{2u-1}}\mu(e)\\
&=&-\sum_{e|m,\,(m,Q(e))=1} 1\\
&\leq & 0\,,
\end{eqnarray*}
since $e\in C_{2u-1}$ implies that $\nu(e)=2u-1$ and hence that
$\mu(e)=-1$.  Thus
$$U=\sum^{\infty}_{u=1}\bg{\twosum{d\in\,B_{2u-1}}{d|m}\mu(d)}\leq 0\,,$$ 
so that
\begin{eqnarray*}
\sum_{d|m}\mu^+(d)&=&\sum_{d|m,\,d\in T^+(y)}\mu(d)\\
&=&\sum_{d|m}\mu(d)-U\\
&\geq&\sum_{d|m}\mu(d)
\end{eqnarray*}
as claimed. The inequality 
$$\sum_{d|m}\mu^{-}(d)\leq \sum_{d|m}\mu(d)$$ 
can be proved in the same way.

This completes the proof of Lemma 5.1.  However a useful alternative
way of viewing the combinatoric facts used in the argument is as
follows.  We have
\begin{eqnarray}
S(\cl{A},\cl{P};z)&=&\sum_{n\in\cl{A}}\twosum{d|\Pi(\cl{P},z)}{d|n}\mu(d)
\nonumber\\
&=&\sum_{n\in\cl{A}}\twosum{d|\Pi(\cl{P},z)}{d|n}\mu^{+}(d)-
\sum_{u=1}^{\infty}\sum_{n\in\cl{A}}\twosum{d|n,\,d|\Pi(\cl{P},z)}
{d\in B_{2u-1}}\mu(d)\nonumber\\
&=&\sum_{d|\Pi(\cl{P},z)}\mu^{+}(d)\#\cl{A}_d\,+\,
\sum_{u=1}^{\infty}\sum_{n\in\cl{A}}\twosum{e|(n,\Pi(\cl{P},z)),
\;e\in C_{2u-1}}
{(n,\Pi(\cl{P},z),Q(e))=1}1\nonumber\\
&=&\sum_{d|\Pi(\cl{P},z)}\mu^{+}(d)\#\cl{A}_d\,+\,
\sum_{u=1}^{\infty}\twosum{e\in C_{2u-1}}{e|\Pi(\cl{P},z)}
S(\cl{A}_{e},\cl{P};p(e)).
\end{eqnarray}

Thus far, all we have said applies to any predicate ${\mathfrak P}$,
and any sieve problem.  We now specialize to a sieve problem of
dimension 1, and examine (5.3) in the particular case
$\cl{A}=\cl{A}^+$, which we expect to be extremal.  Here we have
\[S(\cl{A}^+,\op;z)=
\frac{X}{e^{\gamma}\,\log(x^{1/s })}\,
F(s)\,+\,O_{N}\bg{\frac{x}{(s-1)(\log x)^{2}}}\]
for $1<s\le N$, by Theorem 4.1.  Moreover the estimate (4.2) shows that
\[\sum_{d|\Pi(\op,z)}\mu^{+}(d)\#\cl{A}^+_d=
X\sum_{d|\Pi(\op,z)}\frac{\mu^{+}(d)}{d}+O(x(\log x)^{-2}).\]
We therefore conclude that
\[\sum_{d|\Pi(\op,z)}\frac{\mu^{+}(d)}{d}=
\frac{F(s)}{e^{\gamma}\log z}+
\sum_{u=1}^{\infty}\sum_{e\in C_{2u-1}}S(\cl{A}^+_{e},\op;p(e))
+O_N((s-1)^{-1}(\log x)^{-2}),\]
where $y=E_1(x)$, $z=x^{1/s}$ and $1<s\leq N$.  If we replace $s=(\log
x)/(\log z)$ by $s'=(\log y)/(\log z)$ then the right hand side above
is
\begin{equation}
\frac{F(s')}{e^{\gamma}\log z}+
\sum_{u=1}^{\infty}\twosum{e\in C_{2u-1}}{e|\Pi(\cl{P},z)}
S(\cl{A}^+_{e},\op;p(e))
+O_N((s'-1)^{-1}(\log x)^{-3/2})
\end{equation}
for $1<s'\le N$.  We then re-define $s$ as $(\log y)/(\log z)$.  This
produces an upper bound for the sum
\begin{equation}
\sum_{d|\Pi(\op,z)}\frac{\mu^{+}(d)}{d}
\end{equation}
which involves $F(s)$ together with information about the property
${\mathfrak P}$ incorporated in the definition of the sets $C_r$.
Since our goal is to minimize the sum (5.5), we aim to choose
${\mathfrak P}$ so that $S(\cl{A}^+_{e},\op;p(e))$ is as close to 0 as
possible for $e\in C_{2u-1}$.  However the relevant integers $e$ all
have $\nu(e)=2u-1$, so that $\nu(m)$ is even for any
$m\in\cl{A}^+_e$.  Moreover, every such $m$ satisfies $m\le x/e$.
Hence we would have $S(\cl{A}^+_{e},\op;p(e))=1$ providing that
$x/e<p(e)^2$.  Looking back at the definition of the set $C_{2u-1}$ we
see that we would want to have
\[p_1 p_2\ldots p_{2u-3}p_{2u-1}^3>x\]
whenever ${\mathfrak P}(p_1,\ldots,p_{2u-1})$ is false.  Making a
marginal adjustment to produce a condition which involves $y$ rather
than $x$ we therefore take the property ${\mathfrak P}(p_1,\ldots,p_r)$
to say that 
\[p_1 p_2\ldots p_{r-1}p_r^3<y,\]
whence 
\[S^{+}(y)=\bigcap_{t=1}^{(r+1)/2}
\{d\in\N:\;\mu(d)^{2}=1,\;(d=p_{1}...p_{r}\,\Rightarrow
p_{1}p_{2}...p^{3}_{2t-1}<y)\}\,.\] 

Although we have not made $S(\cl{A}^+_{e},\op;p(e))$ completely vanish
for $e$ in $C_{2u-1}$ it can be shown that this construction does indeed
make the sum in (5.4) suitably small.  We have therefore produced an
admissible set of upper bound sieve coefficients $\mu^+(d)$ which
match up with the Selberg sequence $\cl{A}^+$, and the linear
programming argument then shows that both are optimal.

One can discuss the lower bound problem in exactly the same way, using
the sequence $\cl{A}^-$, and leading to the choice 
\[S^{-}(y)=\bigcap_{t=1}^{r/2}\;
\{d\in\N:\;\mu(d)^{2}=1,\;(d=p_{1}...p_{r}\,\Rightarrow
p_{1}p_{2}...p^{3}_{2t}<y)\}\,.\] 
The construction of $\mu^{\pm}(d)$ we have been led to is known as the
{\it Rosser-Iwaniec} sieve for dimension 1, there being variants in
other dimensions.  (The reader should note that, 
except for dimensions $1$ and $1/2$, the general
Rosser-Iwaniec sieve is not known to be optimal.  Indeed in many case
it is known not to be optimal.)  One noteworthy feature of the
construction is that the definition of the weights $\mu^{\pm}(d)$ does
not involve either the parameter $z$ or the function $\omega(d)$.

Although our discussion has been concerned with the case
$\omega(d)=1$, the Rosser-Iwaniec weights may be applied to the general
sieve problem of dimension 1.  Thus if we set
$$M^{\pm}(z,y)=\sum_{d|\Pi(\cl{P},z)}\frac{\omega(d)\,\mu^{\pm}(d)}{d}$$
we will have
\begin{equation}
S(\cl{A},\cl{P};z)\leq
M^{+}(z,y)X+\sum_{d|\Pi(\cl{P},z)}|R_{d}|
\end{equation}
and
\begin{equation}
S(\cl{A},\cl{P};z)\geq M^{-}(z,y)X-\sum_{d|\Pi(\cl{P},z)}|R_{d}|.
\end{equation}
Iwaniec \cite{LS} has established the following bounds for $M^{+}(z,y)$
and $M^{-}(z,y)$.
\begin{theorem} 
Suppose that
$$\sum_{w\leq p<z}\frac{\omega(p)\log p}{p}\leq \log\frac{z}{w}+O(1)$$ 
and write, as usual,
$$W(z)=\prod_{p<z}\bg{1-\frac{\omega(p)}{p}}\,.$$ 
Then
$$M^{+}(z,y)\leq W(z)\{F(s)+O(e^{-s}(\log y)^{-1/3 })\}$$ 
and
$$M^{-}(z,y)\leq W(z)\{f(s)+O(e^{-s}(\log y)^{-1/3 })\}\,.$$
\end{theorem}
We shall not prove this theorem here.

{\bf Remark.} Note that Iwaniec requires only a one-sided condition,
in contrast to the two-sided condition (3.7) introduced in the context
of the Selberg sieve.  Thus Theorem 5.1 applies to sieves of dimension
less than 1, and even to certain problems without a well-defined dimension.

In view of (5.6) and (5.7), one obtains the following inequalities.
\begin{cor} Under the condition of Theorem 5.1 we have
\[
S(\cl{A},\cl{P};z)\leq  XW(z)\{F(s)+O(e^{-s}(\log
y)^{-1/3 })\}+\sum_{d<y}|R_{d}|
\]
and
\[
S(\cl{A},\cl{P};z)\geq XW(z)\{f(s)+O(e^{-s}(\log y)^{-1/3 })\} -
\sum_{d<y}|R_{d}|.
\]
\end{cor}
{\bf Example 1.} Let $\cl{A}=\{n\in\N:\,n\leq x\,\}\;,
X=x,\,\omega(d)=1$.  Then
$R_{d}=O(1)$ since
$$\#\cl{A}_{d}=\big[\frac{x}{d}\big]=\frac{X\omega(d)}{d}+O(1)\,.$$
From Corollary 5.1 one obtains
\begin{equation}
S(\cl{A},\op;z)\geq
x\prod_{p<z}\bg{1-\frac{1}{p}}(f(s)+o(1))-\sum_{d<y}|R_{d}|.
\end{equation}
Let $y=\frac{x}{(\log x)^{2}}$ and suppose that
$y^{1/4 }<z<y^{1/2}$. Then
$$2<s=\frac{\log y}{\log z} <4,$$
whence
$$f(s)=\frac{2e^{\gamma}\log (s-1)}{s}.$$ 
Thus (5.8), in conjunction with Mertens Theorem (1.4), gives 
\begin{eqnarray*}
S(\cl{A},\op;z)&\geq &x\{\frac{e^{-\gamma}}{\log z}+O(\frac{1}{(\log z)^2})\}
\left(\frac{2e^{\gamma}\log (s-1)}{s}+O((\log y)^{-1/3})\right)\\
&&\hspace{4cm}+O\bg{\frac{x}{(\log x)^{2}}}\\
&=&\frac{2x}{\log y}\{1+o(1)\}\log (s-1)+O\bg{\frac{x}{(\log x)^{2}}}.
\end{eqnarray*}
Hence
$$S(\cl{A},\op;z)\gg \frac{x}{\log x}$$
for any constant value of $s$ strictly greater than 2.
To detect primes, one would need to consider the situation 
with $s=2$. We just fail to find primes, which is not 
surprising since the sequence $\cl{A}^{-}$ contains no primes 
(the {\it parity phenomenon!}). 

{\bf Example 2.} Let $0<\theta<1$, and choose
$$\cl{A}=\{\,n\,:\,x-x^{\theta}<n\leq x\,\},$$ 
$X=x^{\theta},\;\omega(d)=1$, and
$$y=\frac{x^{\theta}}{\log x},\; z=x^{\frac{\theta}{2}-\delta}$$
with $\delta >0$, so that
$$s=\frac{\log y}{\log z}=\frac{2\theta}{\theta - 2\delta}+o(1)>2$$
and $R_{d}=O(1)$. As in Example 1, it follows that
$$S(\cl{A},\op;z)\gg \frac{x^{\theta}}{\log x}\,.$$ 
Hence, if $x$ is large enough, the interval $x-x^{\theta}<n\leq x$
contains at least one integer $n$ all of whose prime factors $p$
satisfy $p\geq x^{\theta/2-\delta}$.  In particular if
$r\in\N$ and $\theta>2/(r+1)$, then we may choose $\delta>0$ so that 
$\theta/2-\delta>1/(r+1)$.  Thus we will have $p>x^{1/(r+1)}$ so that
$n$ can have at most $r$ prime factors.  

In general we say that a positive
integer $n$ is an {\it almost prime} of type $P_r$, if it has at most
$r$ prime factors, counted according to multiplicity.  We may then
conclude that if $\theta>2/(r+1)$, and if $x$ is sufficiently large,
then the interval $x-x^{\theta}<n\leq x$ contains at least one $P_r$
number.  The necessary size for $\theta$ can be reduced (see Example 1
after Theorem 6.2) and
it is an interesting problem to know just small it may be taken.

{\bf Example 3} (The {\it twin primes problem}). Let
$$\cl{A}=\{\,p+2\,:\,p\in\,\op,\;p\leq x\,\},\;$$ 
and take $X=\pi(x)$ and
\[\omega(p)=\left\{\begin{array}{cc} \frac{p}{p-1}, &  p>2,\\
0, & p=2.\end{array}\right.\]
In view 
of the Bombieri-Vinogradov Theorem, one can take 
$y=x^{1/2}(\log x)^{-c_2}$ to obtain
$$\sum_{d<y}|R_{d}|=O\bg{\frac{x}{(\log x)^{3}}}\,.$$ As above, one 
then concludes that
$$S(\cl{A},\op;z)\gg \frac{x}{(\log x)^{2}}\quad \mbox{if}\;s>2\;.$$ 
Since $s=(\log y)/(\log z)$ we may therefore use $z=x^{\theta}$ with 
any constant exponent $\theta< 1/4$.
It follows then that the sequence $\cl{A}$ contains a growing number
of $P_4$ integers as $x$ tends to infinity.
As we shall see in Theorem 6.3, this has been improved by Chen
\cite{Chen} who shows that the same is true for $P_2$ integers.

{\bf Example 4.} Let $\cl{A}=\{\,n^{2}+1\,:\,n\leq x\,\}$, and take
$X=x$ and
$$\omega(p)=\left\{\begin{array}{cc}2,& p\equiv 1\mod{4},\\
0, & p\equiv-1\mod{4},\\
1,&  p=2. \end{array}\right.$$
Then $R_{d}=O(\omega(d))$ and
$$\sum_{w\leq p<z}\frac{\omega(p)\log p}{p}=\log\frac{z}{w}+O(1).$$ 
If we take $y=x(\log x)^{-2}$ it follows that
\[\sum_{d<y}|R_{d}|=O\bg{\frac{x}{(\log x)^{2}}}.\]
Thus if $z=x^{\theta}$ with a constant exponent $\theta>1/2$ we will
find that $s>2$ and hence
$$S(\cl{A},\op;z)\gg \frac{x}{\log x}.$$ 
It follows, on taking $\theta>1/5$, that there are infinitely many
$P_4$ numbers of the form $n^2+1$.

In a similar way one can prove that if $f$ is an irreducible integer
polynomial such that the values $\{\,f(n)\,|\,n\in\,\N\,\}$ contain 
no common factor, then $f(n)=P_{2k}$ infinitely often.
Here also
improvements are possible.
\bigskip

We conclude this section with the following important result. 

\noindent {\bf The ``Fundamental Lemma"} {\it Suppose that
$$\sum_{d<y}|R_{d}|\ll \frac{x}{(\log x)^{2}}\,.$$ 
Then
$$S(\cl{A},\cl{P};z)\sim XW(z),$$
if $\log z=o(\log y)$, that is, if $s\rightarrow\,\infty$.} 

{\bf Proof.} It suffices to observe that $F(s)=1+O(e^{-s})$ and 
$f(s)=1+O(e^{-s})$, as $s\,\rightarrow\,\infty\,.$

{\bf Remarks.} 

1) This result should be compared with Corollary 2.1, in which one
required $z\le\log y$.

2) Since one can often choose $y=X^{\delta}$ for some $\delta > 0$
it follows from the Fundamental Lemma that
$$S(\cl{A},\cl{P};z)\sim XW(z),$$ 
if $z=X^{\ep(X)}$ with $\ep(X)\,\rightarrow 0$ as $X\,\rightarrow\,\infty$.
By considering the 
sequences $\cl{A}^{\pm}$ one sees that one cannot have such a result when
$\ep(X)\,\not\rightarrow 0$.

3) The significance of the Fundamental Lemma is that one can sieve out
``small primes'' (namely those below $X^{\ep(X)}$) as an initial
stage in some more complicated argument, and still have an asymptotic formula. 
In order to make use of this information one usually wants a
quantitative form of the Fundamental Lemma, but this is easily established.

4) One can obtain analogous results on the weaker assumption that
$$\sum_{p\leq z}\frac{\omega(p)\log p}{p}\ll \log z$$ 
as $z\,\rightarrow\,\infty$.

\renewcommand{\thesection}{6}
\section{The weighted sieve}

If the sequence $\cl{A}$ contains only positive integers $n\le N$, and
we can show that $S(\cl{A},\op,N^{\theta})>0$ for some
$\theta>1/(r+1)$, then we can conclude that $\cl{A}$ contains at least
one $P_r$ number.  However one can often derive better results by
using a {\em weighted sieve}, in which certain $P_r$ numbers with
prime factors below $N^{1/(r+1)}$ are also counted.

In general we let $N=\max_{n\in\,\cl{A}}|n|$ and we choose constants
$0<\alpha <\beta$.  We then set
$$W:=W(\cl{A},\cl{P};\alpha,\beta)=
\twosum{n\in\cl{A}}{(n,\Pi(\cl{P},N^{\alpha}))=1}
\bg{1-\twosum{p|n}{N^{\alpha}\leq p<N^{\beta}} w_{p}}\;,$$
where the weights $ w_p\ge 0$ are to be chosen so that
\begin{equation}
\twosum{p|n}{N^{\alpha}\leq p<N^{\beta}} \!\!w_{p}\;\ge\;
1\;\;\;\mbox{if}\;\;\;\nu(n)\ge r+1.
\end{equation}
This condition ensures that
\begin{eqnarray*}
W&\le& \twosum{n\in\cl{A},\,\nu(n)\le r}{(n,\Pi(\cl{P},N^{\alpha}))=1}
\bg{1-\twosum{p|n}{N^{\alpha}\leq p<N^{\beta}} w_{p}}\rule[-2cm]{0cm}{2cm}\\
&\le& \sum_{n\in\cl{A},\,\nu(n)\le r}1.
\end{eqnarray*}
If we know that
\[\#\{n\in\cl{A}:\,\exists p^2|n,\, N^{\alpha}\le p<N^{\beta}\}\ll
\frac{X}{(\log N)^2},\]
we can deduce that
\[W\le\#\{n\in\cl{A}:\,n=P_r\}+O(\frac{X}{(\log N)^2}).\]
Thus if we can also show that
\begin{equation}
W\gg\frac{X}{\log N},
\end{equation}
we will be able to deduce that
\[\#\{n\in\cl{A}:\,n=P_r\}\gg\frac{X}{\log N}.\]

The optimal choice for the weights $ w_p$ is not known, however the
choice
$$ w_{p}=\frac{\beta}{(r+1)\beta-1}\bg{1-\frac{\log p}{\beta\,\log N}}\,$$
leads to some fairly satisfactory results.  These are known as
{\em Richert's logarithmic weights}.  

We shall assume henceforth that
$$(r+1)\beta-1>0.$$
Then if $p<N^{\beta}$ we will have $ w_{p}\ge 0$ as required. Now
suppose that $n\in\cl{A}$ with $\text{h.c.f.}(n,\Pi(\cl{P},N^{\alpha}))=1$, and
consider the sum
$$S=\twosum{p|n}{N^{\alpha}\leq p<N^{\beta}} w_{p}\,.$$  
Then
\begin{eqnarray*}
S&=&\frac{\beta}{(r+1)\beta-1}{\twosum{p|n}{N^{\alpha}\leq
p<N^{\beta}}}\bg{1-\frac{\log p}{\beta\,\log N}}\\
&\geq &\frac{\beta}{(r+1)\beta-1}\sum_{p|n}\bg{1-\frac{\log
p}{\beta\,\log N}},
\end{eqnarray*} 
since $n$ has no prime factors $p<N^{\alpha}$, and
$$1-\frac{\log p}{\beta\,\log N}\leq 0$$
for any prime factor $p\geq N^{\beta}$.  However
\[\sum_{p|n}\log p\le \log |n|\le\log N,\]
so that
\[\sum_{p|n}\bg{1-\frac{\log p}{\beta\,\log N}}\ge
\nu(n)-\frac{1}{\beta}.\]
Hence if $\nu(n)\ge r+1$ we will have
\[S\geq \frac{\beta}{(r+1)\beta-1}\bg{\nu(n)-\frac{1}{\beta}}
\geq 1,\]
as required for (6.1).

We now examine the estimate (6.2).  By definition,
$$W=S(\cl{A},\cl{P};N^{\alpha})-\sum_{N^{\alpha}\leq
p<N^{\beta}} w_{p}\,S(\cl{A}_{p},\cl{P};N^{\alpha})\,.$$
We plan to apply Corollary 5.1.  We therefore make the assumption that
\[\sum_{d<N^{\gamma}}|R_{d}|\ll \frac{X}{(\log N)^{2}}\] 
for some fixed $\gamma>0$, and we set $y=N^{\gamma}$.  
According to Corollary 5.1 we will then have
$$S(\cl{A},\cl{P};N^{\alpha})\geq XW(N^{\alpha})\,
\bg{f\bg{\frac{\gamma}{\alpha}}+o(1)}.$$
Moreover, since
$$(\cl{A}_{p})_{d}=\frac{X\omega(p)}{p}\,\frac{\omega(d)}{d}+R_{pd}, $$
it follows that
$$S(\cl{A}_{p},\cl{P};N^{\alpha})\leq \frac{X\omega(p)}{p}\,W(N^{\alpha})\,
\bg{F\bg{\frac{\log N^{\gamma}/p}{\log
N^{\alpha}}}+o(1)}+\sum_{d<\frac{N^{\gamma}}{p},\;d|\Pi(\cl{P},N^{\alpha})}|R_{pd}|\,.$$
Moreover,
$$\sum_{N^{\alpha}\leq p<N^{\beta}} w_{p}
\sum_{d<N^{\gamma}/p,\,d|\Pi(\cl{P},N^{\alpha})}|R_{pd}|\ll
\sum_{k<N^{\gamma}}|R_{k}|\ll \frac{X}{(\log N)^{2}},$$
since an integer $k<N^{\gamma}$ can have at most $\gamma/\alpha$ prime
factors $p\ge N^{\alpha}$.  It follows that
\begin{eqnarray*}
\lefteqn{\sum_{N^{\alpha}\leq
p<N^{\beta}} w_{p}\,S(\cl{A}_{p},\cl{P};N^{\alpha})}\\
&\leq &XW(N^{\alpha})\sum_{N^{\alpha}\leq
p<N^{\beta}} w_{p}\,\frac{\omega(p)}{p}\,\{F\bg{\frac{\log
N^{\gamma}/p}{\log N^{\alpha}}}+o(1)\} +O(\frac{X}{(\log N)^{2}})\\
&\leq &XW(N^{\alpha})\sum_{N^{\alpha}\leq
p<N^{\beta}} w_{p}\,\frac{\omega(p)}{p}F\bg{\frac{\log
N^{\gamma}/p}{\log N^{\alpha}}} +o(\frac{X}{\log N}).
\end{eqnarray*}
We therefore conclude that
\[W\geq XW(N^{\alpha})\{f\bg{\frac{\gamma}{\alpha}}-
\sum_{N^{\alpha}\leq
p<N^{\beta}} w_{p}\,\frac{\omega(p)}{p}F\bg{\frac{\log
N^{\gamma}/p}{\log N^{\alpha}}}\} +o(\frac{X}{\log N}).\]

We summarize our conclusions as follows.
\begin{theorem}
Suppose that
$$\cl{A}\subseteq \,\Z\cap [-N,N]$$
and that
$$ \#\cl{A}_{d}=X\,
\frac{\omega(d)}{d}+R_{d}\,,$$
and assume that the following conditions hold.
\begin{itemize}

\item[(i)]
$$ \sum_{z<p\leq w}\frac{\omega(p)\log p}{p}=\log w/z+O(1),\;\;\;
(2\leq z\leq w);$$

\item[(ii)] 
$$ \beta>\alpha>0,\;\;\; (r+1)\beta-1>0\,;$$

\item[(iii)] 
\[\#\{n\in\cl{A}:\,\exists p^2|n,\, N^{\alpha}\le p<N^{\beta}\}\ll
\frac{X}{(\log N)^2};\]

\item[(iv)]
$$ \sum_{d<N^{\gamma}}|R_{d}|\ll \frac{X}{(\log N)^{2}};$$ 

\item[(v)] 
$$ f\bg{\frac{\gamma}{\alpha}}-\sum_{N^{\alpha}\leq
p<N^{\beta}}w_p\frac{\omega(p)}{p}\,F\bg{\frac{\log
N^{\gamma}/p}{\log N^{\alpha}}}\gg 1$$
where
$$w_{p}=\frac{\beta}{(r+1)\beta-1}\bg{1-\frac{\log p}{\beta\log
N}}\,.$$ 
\end{itemize}
Then the sequence $\cl{A}$ contains $\gg \frac{X}{\log N}\,,$
numbers of type $P_r$.
\end{theorem}

Our task now is to examine condition (v) in the above theorem.
Let
$$S(t)=\sum_{N^{\alpha}\leq p<t}\frac{\omega(p)\log p}{p}$$
and
$$h(p)=\frac{w_{p}}{\log p}\,F\bg{\frac{\log
N^{\gamma}/p}{\log N^{\alpha}}}\,,$$ 
so that
\[\sum_{N^{\alpha}\leq
p<N^{\beta}}w_p\frac{\omega(p)}{p}\,F\bg{\frac{\log
N^{\gamma}/p}{\log N^{\alpha}}}=
\sum_{N^{\alpha}\leq p<N^{\beta}}\frac{\omega(p)\log p}{p}\,h(p).\]
By partial summation we find that
\begin{eqnarray}
\sum_{N^{\alpha}\leq p<N^{\beta}}\frac{\omega(p)\log p}{p}\,h(p) 
&=&\big[S(t)h(t)\big]^{N^{\beta}}_{N^{\alpha}}- 
\int^{N^{\beta}}_{N^{\alpha}}S(t)h'(t)\,dt\nonumber\\
&=& -\int^{N^{\beta}}_{N^{\alpha}}S(t)h'(t)\,dt,
\end{eqnarray}
since $S(N^{\alpha})=0$ and $h(N^{\beta})=0$.  According to assumption 
(i) we have $S(t)=\log t-\alpha\log N+O(1)$.  The contribution to (6.3)
arising from the error term is
\[\ll\int^{N^{\beta}}_{N^{\alpha}}|h'(t)|\,dt=
\int^{\beta}_{\alpha}|\frac{d\,h(N^v)}{d\,v}|\,dv.\]
However
\[h(N^v)=
\frac{1}{\log N}
\frac{\beta}{(r+1)\beta-1}\big(\frac{1}{v}-\frac{1}{\beta}\big)
F\big(\frac{\gamma-v}{\alpha}\big),\]
whence
\[\frac{d\,h(N^v)}{d\,v}\ll (\log N)^{-1}\]
uniformly for $\alpha\le v\le\beta$.  We therefore conclude that (6.3)
is
\begin{eqnarray*}
\lefteqn{ -\int^{N^{\beta}}_{N^{\alpha}}\{\log t-\alpha\log N+O(1)\}
h'(t)\,dt}\\
&=& -\big[\{\log t-\log N^{\alpha}\}h(t)\big]^{N^{\beta}}
_{N^{\alpha}}+\int^{N^{\alpha}}_{N^{\alpha}}\frac{h(t)}{t}dt\,
+O\big(\frac{1}{\log N}\big)\\
&=&\frac{\beta}{(r+1)\beta-1}\int^{\beta}_{\alpha}
\big(\frac{1}{v}-\frac{1}{\beta}\big)
F\big(\frac{\gamma-v}{\alpha}\big)\,dv+O\big(\frac{1}{\log N}\big).
\end{eqnarray*}
Hence the inequality
\begin{equation}
f\bg{\frac{\gamma}{\alpha}}>\frac{\beta}{(r+1)\beta-1}\int^{\beta}_{\alpha}
\big(\frac{1}{v}-\frac{1}{\beta}\big)
F\big(\frac{\gamma-v}{\alpha}\big)\,dv
\end{equation}
is necessary and sufficient for condition (v) of Theorem 6.1

In order to express $f$ and $F$ in terms of elementary functions we
shall impose the condition $\gamma/4\le\alpha\le\gamma/2$.  For this range we
will have
\[f(\frac{\gamma}{\alpha}) =
2e^{\gamma_{0}}\frac{\log(\gamma/\alpha-1)}{\gamma/\alpha}\]
and
\[F\big(\frac{\gamma-v}{\alpha}\big)=
2e^{\gamma_{0}}\frac{\alpha}{\gamma-v},\]
where we have written $\gamma_0$ for Euler's constant, to avoid
confusion with the parameter $\gamma$. Thus (6.4) is equivalent to the
condition
$$\log\bg{\frac{\gamma}{\alpha}-1}> \frac{\beta}{(r+1)\beta-1}
\int^{\beta}_{\alpha}\frac{\gamma}{\gamma-v}
\bg{\frac{1}{v}-\frac{1}{\beta}}\,dv\,.$$
We can now perform the
integration on the right hand side, to obtain
\begin{equation}
\log \bg{\frac{\gamma}{\alpha}-1}>\frac{1}{(r+1)\beta-1}
\bg{\beta\log \frac{\beta}{\alpha}-(\gamma-\beta)\log
\frac{\gamma-\alpha}{\gamma-\beta}}\,.
\end{equation}
We shall choose
$\alpha=\gamma/4$ and
\[\beta=\frac{\gamma}{1+3^{-r}}.\]
(These are in fact optimal, as a relatively easy calculation shows.
However we do not need to know that the choice is optimal to proceed.)
The above values are compatible with condition (ii) of Theorem 6.1 
providng that
\begin{equation}
\gamma>\frac{1+3^{-r}}{r+1}.
\end{equation}
Moreover (6.5) then reduces to
\[\gamma>\frac{1}{r+1-\frac{\log 4/(1+3^{-r})}{\log 3}},\]
which is a stronger condition than (6.6).

We therefore have the following result.
\begin{theorem}
Suppose the assumptions of Theorem 6 hold, with
\[\alpha=\gamma/4,\;\;\;\mbox{and}\;\;\;\beta=\frac{\gamma}{1+3^{-r}},\]
and with condition (v) replaced by
$$\gamma>\frac{1}{\Lambda_{r}}\;,$$
where
$$\Lambda_{r}:=r+1-\frac{\log 4/(1+3^{-r})}{\log 3}\;.$$
Then the sequence $\cl{A}$ contains $\gg \frac{X}{\log N}\,,$
numbers of type $P_r$.
\end{theorem}

{\bf Remarks.} 

1) For $r\geq 2$ we have $r-\frac{2}{7}<\Lambda_{r}
<r-\frac{1}{7}$.  In particular we have $\Lambda_{2}\geq \frac{11}{6}$.

2) The only parameters which enter into the theorem in a crucial way
are $N$, which measures the size of elements of $\cl{A}$, and $\gamma$
which measures the size of the remainders, in terms of $N$.  The
parameter $\gamma$ is often called the ``level of distribution'' (or
more precisely, since we may not know the optimal value for $\gamma$,
an ``admissible level of distribution'').

{\bf Example 1.} Let
$$\cl{A}=\{n\in\,\N:\;x-x^{\theta}<n\leq x\,\},$$ 
and take $\cl{P}=\op$, $X=x^{\theta}$ and $\omega(p)=1$.  Then 
$R_{d}\ll 1$ so that we may choose any $\gamma<\theta$. 
Then the assumptions of Theorem 6.2 hold true providing that
$\gamma>\Lambda_r^{-1}$. We therefore conclude that $\cl{A}$ contains a
$P_r$ almost-prime if $x$ is large enough, providing that
$\theta>\Lambda_r^{-1}$.   

For $r=2$ much stronger results are known.  According to work of
Baker, Harman and Pintz \cite{BHP}, the sequence $\cl{A}$ actually
contains a prime, for the exponent $\theta=0.525$, which is smaller
than $\Lambda_2^{-1}$.  Moreover Liu \cite{Liu} has show that there
are $P_2$'s as soon as $\theta\ge 0.436$.  It would be nice to know
that $\theta>1/r$ sufficed to ensure the existence of $P_r$'s in
$\cl{A}$, for every $r$.

{\bf Example 2.} Let $N$ be an even integer and put
$$\cl{A}=\{\,N-p\,:\,p\in\,\op,\;3\leq p\leq N-3 \,\},$$
$$\cl{P}=\{p\in\op:\,p\nmid N\},$$
$$\omega(p)=\twos{0,}{\; \mbox{if}\;p|N, }{\frac{p}{p-1},}
{\; \mbox{if}\;p\nmid N,}\,$$
and $X=\li(N)$.  As in our discussion of this example
in \S 3, we find, via the Bombieri-Vinogradov Theorem, that any
$\gamma<1/2$ will be admissable for the remainder sum.  Since
$\Lambda_3>2$ this suffices to show that$\cl{A}$ contains a $P_3$ for
large enough $N$, so that every sufficiently large even integer $2n$
may be written as a sum of a prime and a $P_3$ almost prime.  

In this second example we see that $\Lambda_2$ is only just less 
than 2, so we come quite close to
handling $P_2$'s this way.  However to achieve this requires an
ingenious new idea.

{\bf Chen's theorem.} {\it Every sufficiently large even integer 
$N$ is a sum of a prime and a $P_2$ almost-prime.

More precisely, for every sufficiently large positive 
integer $N$ we have
$$\#\{p\in\,\op:\;N-p\,\in\,\op_{2}\,\}\geq 0.335\,C_{2}\, 
(\prod_{p|N,\, p\neq 2}\frac{p-1}{p-2})\,\frac{N}{(\log N)^{2}},$$
where
$$C_{2}:= 2\prod_{p|N,\, p\neq 2}(1 - \frac{1}{{(p-1)}^{2}})$$
as in} \S\;3.

{\bf Sketch proof.} Given an even positive integer $N$, let
$$\cl{A}=\{\,N-p\,:\,p\in\,\op,\;3\leq p\leq N-3\,\}.$$
Define $\cl{P}$ and $\omega(d)$ as in the previous example,
and let $z>2$.  For $n\in\,\N$, let $a(n)=1$ if 
\begin{equation}
n=p_1p_2p_3,\;\;\;\; p_i\;\mbox{prime},\;\;\;\;
p_1< N^{1/3}\le p_2\le p_3,
\end{equation}
and let $a(n)=0$ otherwise.  We then consider the sum
\[S_0:=S(\cl{A},\cl{P};z)-\frac{1}{2}\sum_{z\leq
p<N^{1/3}}S(\cl{A}_{p},\cl{P};z)-\frac{1}{2}
\sum_{n\in\cl{A},\,(n,\Pi(\cl{P},z))=1}a(n).\]
The first two terms of this may be thought of as giving a weighted
sieve, with constant weights $w_p=1/2$.

Write $S_0^{*}$ for the contribution
to $S_0$ arising from those values of $n$ which are not square-free.
If $z$ is a positive power of $N$ then it is easily shown that
\[S_0^{*}\ll \frac{N}{(\log N)^3},\]
which will be negligible.  (In proving this it is useful to 
note that if $p^2|n$ and $(n,\Pi(\cl{P},z))=1$, then $p\ge z$.)

Let $w(n)$ be the weight attached to $n$ in the expression $S_0$.
Clearly we have $w(n)\le 1$ for every $n$.  We claim that $w(n)\le 0$ 
for any square-free integer $n\in\cl{A}$, unless $n$ is a $P_2$.
Subject to this assertion, we will then have
\[S_0\le \#\{n\in\cl{A}:\, n=P_2\}+O(\frac{N}{(\log N)^3}).\]
To verify the claim take a square-free integer $n\in\cl{A}$ with
$w(n)>0$.  Then we will have $(n,\Pi(\cl{P},z))=1$.  Moreover there
can be at most one prime factor $p|n$ in the range $z\le p<N^{1/3}$,
and clearly any integer $n<N$ can have at most two prime factors $p\ge
N^{1/3}$.  Thus if $n$ is not a $P_2$ almost-prime it must be of the
form (6.7), so that $a(n)=1$.  However it is clear that in this case we
have $w(n)=0$.  This establishes the claim.

The terms $S(\cl{A},\cl{P};z)$ and $S(\cl{A}_{p},\cl{P};z)$ are
estimated from below and above respectively, just as in the standard
weighted sieve.  However it is necessary to choose $z$ somewhat smaller
than before, as $z=N^{1/10}$.  As a result one has to evaluate $f(5)$,
for example, by numerical integration.

However the key new ingredient is the treatment of the sum
\[\sum_{n\in\cl{A},\,(n,\Pi(\cl{P},z))=1)}a(n).\]
Hitherto the only information about $\cl{A}$ that we have used comes
from the estimate
\begin{equation}
\#\cl{A}_{d}=X\frac{\omega(d)}{d}+R_{d}.
\end{equation}
However we now use the precise structure of $\cl{A}$
to re-write the sum above as
\[\#\{p\in\cl{B}\}=S(\cl{B},\op;N^{1/2}\}+O(N^{1/2}),\]
where
$$\cl{B}=\{\,N-p_{1}p_{2}p_{3}\,:\,p_1p_2p_3<N,\;z\leq p_{1}< N^{1/3 }
\leq p_{2}< p_3\,\}.$$
Thus we change our attention to a quite different sequence.
This device has been called the {\it``reversal of r\^{o}les"}, or 
{\it``Chen's twist"}.  (One should note however that although Chen's
application of this idea is arguably the most spectacular, the principle was
independently discovered by Iwaniec, amongst others.)

Although the set $\cl{B}$ looks complicated, it is, in fact,
essentially as 
simple as $\cl{A}$. An analogue of the Bombieri-Vinogradov
Theorem can be establised, showing that if
$$\#\cl{B}_{d}=X\frac{\omega(d)}{d}+R_{d}$$
with a suitable value for $X$, then
\[\sum_{d<N^{\gamma}}|R_{d}|\ll \frac{X}{(\log X)^{3}}\]
for any fixed $ \gamma<\frac{1}{2}$.
 
The standard Selberg upper bound for $S(\cl{B},\op;N^{\frac{1}{2}})$
now allows one to complete the proof of a positive lower bound for $S_0$. 
\bigskip

{\bf Remark}

While the parity phenomenon gives a limitation to the power of sieve
methods which are based purely on the relation (6.8), it is no longer
relevant once one uses additional information.  Thus the reversal of
r\^{o}les trick has the potential to circumvent the parity problem.
\bigskip

{\bf Other applications of the reversal of r\^{o}les trick.} 

1) One can show via the circle method that are infinitely many triples
of distinct primes  $p_{1},p_{2},p_{3}$ which form an arithmetic 
progression, so that $p_2-p_1=p_3-p_2$.  On the other hand 
it is an open problem
whether or not there are infinitely many 4-tuples of distinct primes
in arithmetic progression.

However, one can combine the circle method with the sieve, and use 
the reversal of r\^{o}les trick to give infinitely many 4-tuples
$p_{1},p_{2},p_{3},n$ in arithmetic progression with $n$ a $P_2$
almost-prime, (see Heath-Brown \cite{AP}).

2) In Example 3 of \S 5 we applied the Rosser-Iwaniec lower bound
sieve, together with the  Bombieri-Vinogradov Theorem, to the 
sequence 
\[\cl{A}=\{p+2:\,p\in\op,\,p\le x\}.\]
This was enough to show
that if $\theta<1/4$ then 
\begin{equation}
S(\cl{A},\op,x^{\theta})\gg\frac{x}{(\log x)^2}.
\end{equation}
In particular this shows that $\cl{A}$ contains $P_4$ numbers if we
choose $\theta>1/5$.

However by using the reversal of r\^{o}les trick it is possible to
show the existence of an admissible constant $\theta>1/4$ for which
(6.9) still holds, thereby showing that $\cl{A}$ contains $P_3$ numbers,
without the need for a weighted sieve.

\end{document}